\documentclass{amsart}
\begin{document}
\newcommand{\GL}{{\mathrm{GL}}}
\newcommand{\deform}{{\mathcal D}}
\newcommand{\iso}{{\stackrel{\sim}{\rightarrow}}}
\newcommand{\Sym}{{\mathrm{Sym}}}
\newcommand{\supp}{{\mathrm{supp}\,}}
\newcommand{\abel}{{\mathrm{abel}}}
\newcommand{\discr}{{\mathrm{discr}}}
\newcommand{\trace}{ {\mathrm{trace}} }
\newcommand{\Res}{{\mathrm{Res}}}
\newcommand{\one}{{\mathbf 1}}
\newtheorem{Lemma}[subsubsection]{Lemma}
\newtheorem{Theorem}[subsubsection]{Theorem}
\newtheorem{Proposition}[subsubsection]{Proposition}
\newcommand{\ZZ}{{\mathbb Z}}
\newcommand{\dbold}{{\mathbf d}}
\newcommand{\sbold}{{\mathbf s}}
\newcommand{\vbold}{{\mathbf v}}
\newcommand{\CC}{{\mathbb C}}
\newcommand{\abold}{{\mathbf a}}
\newcommand{\barpartial}{{\overline{\partial}}}
\newcommand{\RR}{{\mathcal R}}
\newcommand{\PP}{{\mathbb P}}
\newcommand{\OO}{{\mathcal O}}
\newcommand{\II}{{\mathcal I}}
\newcommand{\diag}{\mathrm{diag}}
\newcommand{\pluck}{\mbox{pl\"{u}ck}}
\newcommand{\lead}{{\mathrm lead}}
\newcommand{\curve}{C}
\newcommand{\jacobian}{J}
\newcommand{\ord}{\mathrm{ord}}

\title{Abeliants and their
application to an elementary construction of Jacobians}
\author{Greg W. Anderson}
\email{gwanders@math.umn.edu}
\address{School of Mathematics, Univ.\ of Minnesota,
Minneapolis, MN 55455, USA}
\date{This paper has been published. The journal citation is
 Advances in Mathematics \textbf{172}(2002), 169--205.}

\maketitle
\begin{abstract}
The {\em abeliant} is a
polynomial rule which to each
$n$ by $n$ by
$n+2$ array with entries in a commutative ring with unit
associates an $n$ by $n$ matrix with entries in the same ring.
The theory of abeliants,  first introduced in an
earlier paper of the author, is simplified and extended here. Now let
$\jacobian$ be the Jacobian of a nonsingular projective algebraic curve
defined over an algebraically closed field. With the aid of the
theory of abeliants we obtain explicit defining equations for $\jacobian$
and its group law.
\end{abstract}

\section{Introduction}

Let $\curve$ be a nonsingular projective algebraic curve defined over 
an algebraically closed field $k$ and let $\jacobian$
be the Jacobian of
$\curve$.  The point of the paper is to give an elementary
construction of $\jacobian$, i.~e., to obtain by  purely algebraic and
relatively simple means explicit defining equations for
$\jacobian$ and its group law. For historical perspective see
\cite{Milne}.
Our construction is similar in spirit to that of
\cite{MumfordTataII}, but differs from the latter in (at least) two
important respects. Firstly,  we need not assume that
$\curve$ is hyperelliptic. Secondly, we obtain a
description of
$\jacobian$ not as a
glued-together collection of affine varieties, but rather
as a projective variety.

Our construction of $\jacobian$ is based in large part on the notion of
{\em abeliant} introduced in the author's earlier paper \cite{Anderson}.
The abeliant is just a
polynomial rule which to each
$n$ by $n$ by
$n+2$ array with entries in a commutative ring with unit
associates an $n$ by $n$ matrix with entries in the same ring. We simplify and extend the
theory of abeliants in this paper (\S\ref{section:Abeliants}). One of the
new results obtained here is an expansion of each entry of the abeliant
 as a sum indexed by four
permutations of $n$ letters (\S\ref{subsection:FourPermutationExpansion}).

Our construction of $\jacobian$ proceeds in three main stages. The first
stage is to derive from the theory of abeliants a theory of {\em
abstract Abel maps} (\S\ref{section:SegreClassification}). The first
stage is more or less pure multilinear algebra and has {\em a priori}
nothing to do with algebraic curves. The abstract Abel map is roughly
analogous to the Pl\"{u}cker embedding, but with this  important
difference: 
 it collapses not $\GL_n$-orbits but rather 
$\GL_n\times
\GL_n$-orbits to lines.
In the second stage of the construction of $\jacobian$ we set up a
representation of divisor classes of
$\curve$ of sufficiently high degree by square rank one matrices with
entries in the function field of
$\curve$ (\S\ref{subsection:Dictionary})
and then set up corresponding matrix representations of addition and
subtraction of divisor classes
(\S\ref{subsection:AdditionAndSubtraction}). Since the matrix
representing a divisor class is well-defined only up to 
$\GL_n\times
\GL_n$-equivalence, we have to solve a problem in the
invariant theory of
$\GL_n\times\GL_n$ to complete the construction of $\jacobian$. Of course
it is precisely this sort of problem that the theory of abstract Abel
maps is designed to solve. Thus the
third and final stage of the construction of $\jacobian$ comes down to
interpreting the abstract Abel map in certain special cases associated
to
$\curve$ as {\em the} Abel map (Theorem~\ref{Theorem:Main}).

The explicit elementary
point of view on hyperelliptic Jacobians developed by Mumford and many
others has been quite useful in number theory and computer
science. To give just two examples of applications, we cite the papers
\cite{FlynnPoonenSchaefer} and \cite{AdlemanDeMarraisHuang}. 
The \linebreak first is a study of the rational points on a certain
curve of genus $2$ connected with
iteration of quadratic polynomials; the second is a cryptologically
motivated study of the discrete logarithm problem in Jacobians of
hyperelliptic curves defined over finite fields.
We expect the explicit elementary point of view on
not-necessarily-hyperelliptic Jacobians developed here to be analogously
useful.  A particularly interesting problem that might be approachable
from our point of view is that of implementing the algorithm of
\cite[Theorem D]{Pila} for finding
$\ell^{th}$ roots of unity modulo $p$; heretofore the sticking
point has  been the lack of a sufficiently explicit model for the
Jacobian of the Fermat curve of degree $\ell$. 

\section{Abeliants}\label{section:Abeliants}
We review, simplify and refine the
theory of {\em abeliants} introduced in the
author's previous paper \cite{Anderson}. 
{\em Rings} are
commutative with unit. 

\subsection{Definition}
\label{subsection:AbeliantDefinition}
 Given an $n$ by $n$ matrix $X$ with entries in some
ring, let
$X^\star$ denote the transpose of the matrix of cofactors of $X$, i.~e.,
the
$n$ by
$n$ matrix  with entry in position $ji$
equal to
$(-1)^{i+j}$ times the determinant of the matrix obtained by
striking row $i$ and column $j$ from
$X$; we then have 
$$X^\star X=XX^\star=\diag\left(\underbrace{\det
X,\dots,\det X}_n\right).$$ Here and elsewhere
$\diag(x_1,\dots,x_n)$ denotes the $n$ by $n$ diagonal
matrix with diagonal entries
$x_1,\dots,x_n$. Now let
$$\left\{X^{(\ell)}\right\}_{\ell=0}^{n+1}$$
be a family of $n$ by $n$ matrices with entries in a
ring $R$ and let
$$\{s_i\}_{i=1}^n\cup \{t_j\}_{j=1}^n$$
be a family  of independent variables.
Following
\cite[p.~496]{Anderson}, we define the {\em abeliant}
$$\abel\left(X^{(0)},\dots,X^{(n+1)}\right)
\;\;\;\left(\mbox{abbreviated
notation:}\;\abel_{\ell=0}^{n+1}X^{(\ell)}\right)$$ of the given
family  of matrices to be
the
$n$ by
$n$ matrix the entry of which in 
position
$ij$ is  the coefficient with which the
monomial
\begin{equation}\label
{equation:AbeliantDefiningMonomial}s_i^{-1}t_j^{-1}\cdot\prod_{a=1}^n
s_a\cdot \prod_{b=1}^n t_b=s_1\cdots \widehat{s_i}\cdots
s_n t_1\cdots \widehat{t_j}\cdots t_n
\end{equation}
appears in the
expansion of the expression
\begin{equation}\label
{equation:AbeliantDefiningExpression}
\trace\left(X^{(0)}\left(\sum_{b=1}^nt_bX^{(b)}\right)^\star
X^{(n+1)}\left(\sum_{a=1}^n s_aX^{(a)}\right)^\star\right)
\end{equation}
as an $R$-linear combination of monomials in the $s$'s and $t$'s.

\subsection{Basic properties}
Let $\left\{X^{(\ell)}\right\}_{\ell=0}^{n+1}$ be a family
of
$n$ by $n$ matrices with entries in a ring $R$. For any
square matrices $X$ and $Y$ with entries in a ring we have 
$$\trace
(XY)=\trace(YX),\;\;\;(XY)^\star=Y^\star X^\star.$$
It follows after a short calculation that
\begin{equation}\label{equation:FundamentalLaw}
\abel_{\ell=0}^{n+1}\left(UX^{(\ell)}V\right) =(\det
U)^2(\det V)^2\abel_{\ell=0}^{n+1}X^{(\ell)}
\end{equation}
for all $n$ by $n$ matrices $U$  and $V$ 
with entries in $R$.  For any  square matrix
$X$ with entries in a ring we have 
$$\trace(X^T)=\trace(X),\;\;\;(X^\star)^T=(X^T)^\star$$ where $X^T$
denotes the transpose of $X$. It follows after a short
calculation that 
\begin{equation}\label{equation:AbeliantTransposeSymmetry}
\abel_{\ell=0}^{n+1}X^{(\langle
0\mid
n+1\rangle
\ell)}=\left(\abel_{\ell=0}^{n+1}X^{(\ell)}\right)^T
=\abel_{\ell=0}^{n+1}\left(X^{(\ell)}\right)^T
\end{equation}
where $\langle 0\mid
n+1\rangle$ denotes the permutation of $\{0,\dots,n+1\}$
exchanging $0$ and $n+1$ and fixing all other
elements. We claim
that
\begin{equation}\label{equation:AbeliantSymmetry}
\left(\abel_{\ell=0}^{n+1}X^{(\pi
\ell)}\right)_{ij} =\left(\abel_{\ell=0}^{n+1}X^{(\ell)}
\right)_{\pi i, \pi j}
\end{equation}
where $\pi$ is any permutation of $\{0,\dots,n+1\}$ fixing $0$ and $n+1$.
We claim further that for $n\geq 2$ we have
\begin{equation}\label{equation:AbeliantDegeneration}
\left(\abel_{\ell=0}^{n+1}X^{([1\mapsto 2]\ell)}
\right)_{12}
=\left(\abel_{\ell=0}^{n+1}X^{(\ell)}\right)_{11}
\end{equation}
where $[1\mapsto 2]$ denotes the mapping of
$\{0,\dots,n+1\}$ to itself sending $1$ to $2$ and
fixing all other elements. Since the proofs of
(\ref{equation:AbeliantSymmetry}) and 
(\ref{equation:AbeliantDegeneration}) are similar, we supply a proof only
for (\ref{equation:AbeliantDegeneration}). To abbreviate notation  let
monomial (\ref{equation:AbeliantDefiningMonomial}) be denoted by
$m_{ij}$, let expression
(\ref{equation:AbeliantDefiningExpression}) be denoted by
$F$, and let 
$$s^\mu  t^\nu=\prod_{a=1}^n s_a^{\mu_a}\cdot \prod_{b=1}^n
t_b^{\nu_b}$$
be a monomial in the $s$'s and $t$'s.  Further, let
$[1\mapsto 2]^*$ be the unique
$R$-algebra endomorphism of the polynomial ring
$R[s_1,\dots,s_n,t_1,\dots,t_n]$ such that 
$$[1\mapsto 2]^*s_a=\left\{\begin{array}{cl}
0&\mbox{if $a=1$,}\\
s_1+s_2&\mbox{if $a=2$,}\\
s_a&\mbox{if $a\geq 3$,}
\end{array}\right.\;\;\;
[1\mapsto 2]^*t_b=\left\{\begin{array}{cl}
0&\mbox{if $b=1$,}\\
t_1+t_2&\mbox{if $b=2$,}\\
t_b&\mbox{if $b\geq 3$,}
\end{array}\right.
$$
for $a,b=1,\dots,n$. 
Now the
coefficient with which the monomial
$m_{12}$ (resp.~$m_{11}$)
appears in the expansion of $[1\mapsto 2]^* F$
(resp.~$F$) as an $R$-linear combination of monomials in the $s$'s
and $t$'s admits interpretation as the left (resp.~right) side of
(\ref{equation:AbeliantDegeneration}).  But the
coefficients in question are equal because 
$$[1\mapsto 2]^*s^\mu t^\nu=
\left\{\begin{array}{ll}
m_{11}+m_{12}+m_{21}+m_{22}&\mbox{if $s^\mu t^\nu=m_{11}$,}\\
\mbox{a polynomial in which $m_{12}$ does not appear}\;\;
&\mbox{if $s^\mu t^\nu\neq m_{11}$.}
\end{array}\right.$$ 
Thus 
claim (\ref{equation:AbeliantDegeneration}) is proved.

\subsection{Evaluation in a special
case}\label{subsection:SpecialEvaluation}
Let $X$, $L$, $M$, and $Q$ be $n$ by $n$ matrices with
entries in a ring
$R$, where $L$,
$M$ and $Q$ are diagonal. We write
$Q=\diag(q_1,\dots,q_n)$.  Let $e^{(\ell)}$
(resp.~$f^{(\ell)}$) be the
$\ell^{th}$ column (resp.~row) of the $n$ by $n$ identity matrix
and put
$$E=\left(e^{(1)}+\cdots+e^{(n)}\right)
\left(f^{(1)}+\cdots+f^{(n)}\right)
=\left[\begin{array}{ccc}
1&\cdots&1\\
\vdots&&\vdots\\
1&\cdots&1\end{array}\right].
$$
We claim that
\begin{equation}\label{equation:SpecialCase}
\abel(X,q_1e^{(1)}f^{(1)},\dots,q_n
e^{(n)}f^{(n)},LEM)=MQ^\star XQ^\star L.
\end{equation}
For the proof we write
$$
S=\sum_{i=1}^n s_i e^{(i)}f^{(i)}=\diag(s_1,\dots,s_n),\;\;\;
T=\sum_{j=1}^n t_j e^{(j)}f^{(j)}=\diag(t_1,
\dots,t_n)$$
where, as above, $s_1,\dots,s_n,t_1,\dots,t_n$ are independent variables.
The identity
$$\begin{array}{cl}
&\displaystyle\trace\left(X\left(\sum_{b=1}^n
q_bt_be^{(b)}f^{(b)}\right)^\star LEM\left(\sum_{a=1}^n
q_as_ae^{(a)}f^{(a)}\right)^\star\right)\\\\
=&\trace(S^{\star}M Q^\star  X
Q^\star  LT^{\star}E)
\end{array}
$$  
is easily verified and suffices to prove the claim.

\subsection{Discriminants}
Let
$\left\{X^{(\ell)}\right\}_{\ell=1}^{n+1}$ be a family
of $n$ by $n$ matrices with entries in a ring $R$. 
 We define the 
{\em discriminant}
$$\Delta\left(X^{(1)},\dots,X^{(n+1)}\right)\;\;
\left(\mbox{abbreviated
notation:}\;\Delta_{\ell=1}^{n+1}X^{(\ell)}\right)$$
of the given family of matrices  to be
the product
$$\left(\det\left(\sum_{\ell=1}^nX^{(\ell)}\right)\right)^{2n-2}\cdot
\prod_{i=1}^n
\left(\det\left(\sum_{\ell\in \{1,\dots,n+1\}\setminus\{i\}}
X^{(\ell)}\right)\right)^2.
$$   
As becomes clear presently, this nonstandard notion of discriminant 
is closely allied with the notion of abeliant.
We have
\begin{equation}
\label{equation:DiscriminantTransformationLaw}
\Delta_{\ell=1}^{n+1}UX^{(\ell)}V=\left(\det
U\right)^{4n-2}\left(\det V \right)^{4n-2}
\Delta_{\ell=1}^{n+1}X^{(\ell)}
\end{equation}
for all $n$ by $n$ matrices $U$ and $V$ with entries
in $R$. We have
\begin{equation}\label{equation:DiscriminantTransposeSymmetry}
\Delta_{\ell=1}^{n+1}X^{(\ell)}=
\Delta_{\ell=1}^{n+1}(X^{(\ell)})^T.
\end{equation}
If 
there exist factorizations
$$X^{(\ell)}=u^{(\ell)}v^{(\ell)}\;\;\;\mbox{for $\ell=1,\dots, n+1$}$$
where $u^{(\ell)}$ (resp.~$v^{(\ell)}$)
is a column (resp.~row) vector with entries in $R$, then we have
$$
\det\left(\sum_{\ell\in
\{1,\dots,n+1\}\setminus\{i\}}X^{(\ell)}\right)=\left|\begin{array}{ccccc}
u^{(1)}&\cdots&\widehat{u^{(i)}}&\cdots&u^{(n+1)}
\end{array}\right|\cdot
\left|
\begin{array}{c}
v^{(1)}\\
\vdots\\
\widehat{v^{(i)}}\\
\vdots\\
v^{(n+1)}
\end{array}\right|
$$
for $i=1,\dots, n+1$ and hence,
given matrices
$$L,M,Q=\diag(q_1,\dots,q_n),e^{(\ell)},f^{(\ell)},E$$
with entries in $R$ as in (\ref{equation:SpecialCase}) above,
we have 
\begin{equation}
\label{equation:DiscriminantSpecialCase}
\Delta\left(q_1e^{(1)}f^{(1)},\dots,q_n e^{(n)}f^{(n)},
LEM\right)=(\det Q)^{4n-4}(\det L)^2(\det M)^2
\end{equation}
after a straightforward calculation we
can safely omit. 

\subsection{The key
relations}\label{subsection:KeyRelation} Let
$\left\{X^{(\ell)}\right\}_{\ell=0}^{n+1}$ be a family
of $n$ by
$n$ matrices with entries in a ring $R$. Assume that
we have factorizations
$$X^{(\ell)}=u^{(\ell)}v^{(\ell)}\;\;\;\mbox{for $\ell=1,\dots,n+1$}
$$
where
$u^{(\ell)}$ (resp.~$v^{(\ell)}$) is a column 
(resp.~row) vector with entries in $R$. 
For the moment we do not
assume that $X^{(0)}$ has such a factorization. We claim
that
\begin{equation}\label{equation:KeyRelation}
\abel_{\ell=0}^{n+1} X^{(\ell)}=
MU^\star X^{(0)}V^\star L,\;\;\;
\Delta_{\ell=1}^{n+1}X^{(\ell)}=\det(MU^\star)^2\det(V^\star
L)^2
\end{equation}
where
$$M=\diag\left(
\left(v^{(n+1)}V^\star\right)_1,\dots,
\left(v^{(n+1)}V^\star\right)_n\right),\;\;\;
U=\left[\begin{array}{ccc}
u^{(1)}&\cdots&u^{(n)}\end{array}\right],
$$\\
$$
V=\left[\begin{array}{c} v^{(1)}\\\vdots\\
v^{(n)}\end{array}\right],\;\;\;
L=\diag\left(
\left(U^\star u^{(n+1)}\right)_1,\dots,\left(U^\star
u^{(n+1)}\right)_n\right).
$$\\
The relations (\ref{equation:KeyRelation})
are key to all applications of the abeliant. Let
$u$ (resp.~$v$) be any column (resp.~row) vector of
length
$n$ with entries in
$R$. We have
\begin{equation}\label{equation:PreHelper}
(U^\star
u)_i=(-1)^{i+1}\left|\begin{array}{cccccc} u
&u^{(1)}&\cdots&\widehat{u^{(i)}}&\cdots&u^{(n)}
\end{array}\right|,\;\;\;
(vV^\star)_j=(-1)^{1+j}\left|\begin{array}{c}
v\\ 
v^{(1)}\\
\vdots\\
\widehat{v^{(j)}}\\
\vdots\\
v^{(n)}
\end{array}\right|
\end{equation}
by Cramer's Rule
and hence
\begin{equation}\label{equation:Helper}
(U^\star u)_i(vV^\star)_i=\det\left(uv+\sum_{\ell\in
\{1,\dots,n\}\setminus\{i\}}X^{(\ell)}\right).
\end{equation}
In particular, we have
$$\det U\cdot\det
V=\det\left(\sum_{\ell=1}^n X^{(\ell)}\right),\;\;\;
(LM)_{ii}=\det\left(\sum_{\ell\in
\{1,\dots,n+1\}\setminus\{i\}}X^{(\ell)}\right).
$$ Repeated application
of (\ref{equation:Helper})  proves the second part of
(\ref{equation:KeyRelation}). Now 
let
$e^{(\ell)}$, $f^{(\ell)}$ and $E$ be as in
(\ref{equation:SpecialCase}) above and put
$$D=\det U\cdot \det V.$$ We have
$$
\begin{array}{cl}
&D^{2n-2}\abel(X^{(0)},u^{(1)}v^{(1)},
\dots,u^{(n)}v^{(n)},u^{(n+1)}v^{(n+1)})\\
=&(\det
U^\star)^2(\det
V^\star)^2\abel(X^{(0)},u^{(1)}v^{(1)},\dots,u^{(n)}v^{(n)},
u^{(n+1)}v^{(n+1)})\\ =&\abel(U^\star X^{(0)}V^\star
,De^{(1)}f^{(1)},\dots,D e^{(n)}f^{(n)}, LEM)\\ =&D^{2n-2}MU^\star
X^{(0)}V^\star L
\end{array}
$$
 by transformation law (\ref{equation:FundamentalLaw}) and
special case (\ref{equation:SpecialCase}). 
The preceding calculation proves the
first part of (\ref{equation:KeyRelation}) provided
that cancellation of the factor $D^{2n-2}$ can be justified. But
there is no loss of generality in assuming that the entries of the
matrix
$X^{(0)}$ and of the vectors
$u^{(\ell)}$
and $v^{(\ell)}$ together constitute a family of
independent variables, and that $R$ is the ring generated by
these variables over the integers; then 
$R$ is an integral domain, $D\neq 0$, cancellation of
the factor $D^{2n-2}$ is permitted, and the
first part of (\ref{equation:KeyRelation}) is proved. 
The proof of  (\ref{equation:KeyRelation}) is complete.

An amplifying remark is now in order. If there exist factorizations
$$X^{(\ell)}=u^{(\ell)}v^{(\ell)}\;\;\;\mbox{for $\ell=0,\dots,n+1$}$$
with
$u^{(\ell)}$ (resp.~$v^{(\ell)}$) a column (resp.~row) vector
with entries in $R$ (notice that $\ell=0$ is now included) then the first
part of  key relation (\ref{equation:KeyRelation}) takes the form
 \begin{equation}\label{equation:KeyRelationCompanion}
\begin{array}{cl}
&\left(\abel_{\ell=0}^{n+1}X^{(\ell)}\right)_{ij}\\\\\\
=&\left|\begin{array}{c}
\widehat{v^{(0)}}\\
\vdots\\
\widehat{v^{(i)}}\\
\vdots\\
\end{array}\right| \cdot
\left|\begin{array}{cccccc}
\cdots&\widehat{u^{(i)}}&\cdots&\widehat{u^{(n+1)}}
\end{array}\right|\cdot
\left|\begin{array}{c}\vdots\\
\widehat{v^{(j)}}\\
\vdots\\
\widehat{v^{(n+1)}}
\end{array}\right|\cdot
\left|\begin{array}{cccccc}
\widehat{u^{(0)}}&\cdots&\widehat{u^{(j)}}&\cdots&
\end{array}\right|
\end{array}
\end{equation}
by (\ref{equation:PreHelper}). Identity
(\ref{equation:KeyRelationCompanion}) is the
method we almost always use for evaluating abeliants of
algebro-geometric interest.

\subsection{Abeliants of
matrices of rank $\leq 1$}
\label{subsection:RankLeqOne}  We say that a matrix $X$ with
entries in some ring is {\em of rank
$\leq 1$} if every two by two submatrix has vanishing
determinant. Now let
$\left\{X^{(\ell)}\right\}_{\ell=0}^{n+1}$ be a
family  of
$n$ by $n$ matrices with entries in a ring $R$.
Assume that $n\geq 2$
and that 
$$\mbox{$X^{(\ell)}$ is of rank $\leq 1$ for $\ell=0,\dots,n+1$.}$$
For distinct
$i,j\in
\{0,\dots,n+1\}$ put
$$D_{ij}=\left(\sum_{\ell\in
\{0,\dots,n+1\}\setminus\{i,j\}}X^{(\ell)}\right).
$$
For $a\in\{0,\dots,n+1\}$, let $[0\mapsto a]$ denote the mapping of
$\{0,\dots,n+1\}$ to itself sending $0$ to $a$ and fixing all other
elements. For distinct $a,b\in\{0,\dots,n+1\}$ let $\langle a\mid
b\rangle$ denote the permutation of $\{0,\dots,n+1\}$
exchanging $a$ and $b$ and fixing all other elements. We
make the following claims:
\begin{equation}\label{equation:RankOneInRankOneOut}
\mbox{$\abel_{\ell=0}^{n+1}X^{(\ell)}$ is a matrix of
rank
$\leq 1$.}
\end{equation}
\begin{equation}
\label{equation:AbeliantDuplicateVanishing}
\left(\abel_{\ell=0}^{n+1}
X^{([0\mapsto a]\ell)}\right)_{ij}=0\;\;\;
\mbox{unless $a\in \{0,n+1\}$ or $i=j=a$.}
\end{equation}
\begin{equation}
\label{equation:RankOneDiagonalAbeliantEvaluationPrime}
\left(\abel_{\ell=0}^{n+1}X^{([0\mapsto n+1]\ell)}\right)_{ij}=
D_{0i}D_{0j}
\end{equation}
\begin{equation}
\label{equation:RankOneDiagonalAbeliantEvaluationBis}
\left(\abel_{\ell=0}^{n+1}
X^{([0\mapsto a]\ell)}\right)_{aa}=D_{0a}D_{0,n+1}
\end{equation}
\begin{equation}
\label{equation:RankOneDiagonalAbeliantEvaluation}
\left(\abel_{\ell=0}^{n+1}
X^{(\ell)}\right)_{aa}=D_{0a}D_{a,n+1}
\end{equation}
\begin{equation}
\label{equation:DiscriminantDefinitionReformulation}
\begin{array}{rccl}
\Delta_{\ell=1}^{n+1} 
X^{(\ell)}&=&&
 \left(\abel_{\ell=0}^{n+1}
X^{(\langle
2|n+1\rangle\langle 0|1\rangle\ell)}\right)_{11}\\\\
&&\cdot &
 \left(\abel_{\ell=0}^{n+1}
X^{(\langle
0|1\rangle\ell)}\right)_{11}\cdot\left(\abel_{\ell=0}^{n+1}X^{(\langle
0|2\rangle\ell)}\right)_{22}\\\\ &&\cdot&
\displaystyle\prod_{a=3}^n\left(\abel_{\ell=0}^{n+1}X^{(\langle
0|a\rangle\ell)}\right)_{aa}^2
\end{array}
\end{equation}
Let
$\left\{\tilde{X}^{(\ell)}\right\}_{\ell=0}^{n+1}$ be
a family of matrices the entries of which constitute a
family of
$(n+2)\cdot n \cdot n$ independent variables. 
Without loss of generality we may
assume that
$R$ is the quotient of the ring generated by the
entries of the
$\tilde{X}^{(\ell)}$ over the integers by the ideal
$I$ generated  by the determinants of all two by two submatrices
of the $\tilde{X}^{(\ell)}$, and we may assume that
$X^{(\ell)}\equiv \tilde{X}^{(\ell)}\bmod{I}$ for all
indices $\ell$. As is well known
the ideal
$I$ is prime and hence the ring $R$ is an integral
domain.  Over  the fraction field
of
$R$ we have factorizations
$X^{(\ell)}=u^{(\ell)}v^{(\ell)}$ with
$u^{(\ell)}$ (resp.~$v^{(\ell)}$) a column
(resp.~row) vector.  The first five claims now follow immediately from
relation (\ref{equation:KeyRelationCompanion}). The last 
claim follows from the penultimate one after a
straightforward calculation we can safely omit. Thus all claims are
proved. 

\subsection{Iterated abeliants}
Let $\left\{X^{(\ell)}\right\}_{\ell=-n-1}^{n+1}$ be a
family of
$n$ by $n$ matrices with entries in a ring $R$.
Assume that all the matrices $X^{(\ell)}$
are of rank $\leq 1$. We claim that
\begin{equation}\label{equation:IteratedAbeliants}
\abel_{\ell=0}^{n+1}
\abel\left(X^{(-\ell)},X^{(1)},\dots,X^{(n+1)}
\right)=\Delta_{\ell=1}^{n+1}X^{(\ell)}\cdot
\abel_{\ell=0}^{n+1}X^{(-\ell)}.
\end{equation}
Arguing in much the same fashion as in the proof of fact
(\ref{equation:RankOneInRankOneOut}) and its companions, we may assume
without loss of generality that  there exist factorizations 
$X^{(\ell)}=u^{(\ell)}v^{(\ell)}$ over $R$ 
with $u^{(\ell)}$ (resp.~$v^{(\ell)}$) a column 
(resp.~row) vector. Then the claim
follows by transformation law
(\ref{equation:FundamentalLaw})
and key relation (\ref{equation:KeyRelation}).

\subsection{Expansion
of the
abeliant}\label{subsection:FourPermutationExpansion} 
Let
$\left\{X^{(\ell)}\right\}_{\ell=0}^{n+1}$ be a family of $n$
by $n$ matrices with entries in a ring $R$. By definition we have
\begin{equation}\label{equation:ExpansionIntermediate}
\left(\abel_{\ell=0}^{n+1}X^{(\ell)}\right)_{ij}=
\sum_{e,f,g,h=1}^n
X^{(0)}_{fg}Y^{(j)}_{gh}X^{(n+1)}_{he}
Y^{(i)}_{ef}
\end{equation}
where
$Y^{(i)}_{ef}$ denotes the
coefficient with which the monomial
$s_1\cdots
\widehat{s_i}\cdots s_n$ appears in the expansion of the
matrix entry
$$\left(\sum_{a=1}^n s_aX^{(a)}\right)^\star_{ef}$$
as an $R$-linear combination of monomials in the $s$'s. Now for any
$n$ by
$n$ matrix
$X$ with entries in a ring we have
$$X^\star_{ef}=\sum_{\pi e=f}(-1)^\pi
\prod_{c\neq e}X_{\pi c,c}$$ where the sum is
extended over  permutations
$\pi$ of $\{1,\dots,n\}$
such that $\pi e=f$,
the product is extended over
$c\in \{1,\dots,n\}\setminus\{e\}$
and $(-1)^\pi$ is the sign of $\pi$. It follows that
$$Y^{(i)}_{ef}=\sum_{\begin{array}{c}
\scriptstyle\pi e=f\\
\scriptstyle\theta e=i\end{array}}(-1)^\pi \prod_{c\neq e}
 X^{(\theta c)}_{\pi c,c}$$
where the sum is extended over permutations $\pi$ and $\theta$
of
$\{1,\dots,n\}$ such that $\pi e=f$ and
$\theta e=i$ and the product is extended over  $c\in
\{1,\dots,n\}\setminus\{e\}$. Substituting
$\theta=\psi^{-1}$, $\pi=\tau\psi^{-1}$, $c=\psi a$ and
making the simplification
$(-1)^{\tau\psi^{-1}}=(-1)^{\tau\psi}$,  we obtain the
expansion
$$Y^{(i)}_{ef}=\sum_{\begin{array}{c}
\scriptstyle\psi i=e\\
\scriptstyle\tau i=f\end{array}}(-1)^{\tau \psi}
\prod_{a\neq 0,i,n+1}X^{(a)}_{\tau
a,\psi a}$$ where the sum is extended over  permutations
$\tau,\psi$ of $\{1,\dots,n\}$
such that $\psi i=e$ and $\tau i=f$ and the product is
extended over 
$a\in \{1,\dots,n\}\setminus\{i\}$. 
Finally, after substituting into
(\ref{equation:ExpansionIntermediate}), we 
obtain the expansion
\begin{equation}
\label{equation:FourPermutationExpansion}
\begin{array}{cl}
&\displaystyle\left(\abel_{\ell=0}^{n+1}X^{(\ell)}
\right)_{ij}\\\\
=&\displaystyle
\sum_{\sigma,\phi,\tau,\psi}
(-1)^{\sigma\phi\tau\psi}
X_{\tau i,\phi j}^{(0)}\cdot
\prod_{b\neq 0,j,n+1}
X^{(b)}_{\sigma b,\phi b}\cdot
X_{\sigma j,\psi i}^{(n+1)}\cdot \prod_{a\neq 0,i,n+1}
X^{(a)}_{\tau a,\psi a}
\end{array}
\end{equation}
 where the sum is extended over   permutations
$\sigma,\phi,\tau,\psi$ of $\{1,\dots,n\}$
and the products are extended over
$a\in \{1,\dots,n\}\setminus\{i\}$ and 
$b\in \{1,\dots,n\}\setminus\{j\}$.

\section{The abstract Abel map}
\label{section:SegreClassification}
We abstract and refine the part of the theory of
\cite{Anderson}  having to do with invariants of $\GL_n\times \GL_n$.

\subsection{Segre matrices}
\label{subsection:SegreSetting}

\subsubsection{Basic data}
Throughout \S\ref{section:SegreClassification} we
work with data
$$k,n,A,L$$ consisting of 
\begin{itemize}
\item an algebraically closed
field
$k$, 
\item an integer $n\geq 2$,
\item a finitely generated
$k$-algebra $A$ without zero divisors, and
\item a finite-dimensional $k$-subspace $L\subset A$.
\end{itemize}
Here and below {\em $k$-algebras} are
commutative with unit. 
We often refer to elements of
$k$ as {\em constants} or {\em scalars}.

\subsubsection{Matrix terminology}
\label{subsubsection:MatrixTerminology}
Let $X$ and $Y$ be matrices with entries in a $k$-algebra $R$.  We say
that
$X$ is \linebreak {\em
$k$-general} if there exists both a row of $X$
with $k$-linearly independent entries
and a column of $X$ with $k$-linearly independent entries.
As in \S\ref{section:Abeliants}, we
say that $X$ is of {\em rank $\leq 1$} if every two by two
submatrix has vanishing determinant.  We say that
$X$ and
$Y$ are {\em
$k$-equivalent} if there exist square matrices $\Phi$ and
$\Psi$ with entries in $k$ such that $\det \Phi\neq 0$, $\det \Psi\neq
0$, the product $\Phi X\Psi$ is defined,  and $Y=\Phi X\Psi$. 
Also, given vectors $x$ and $y$ in a common vector
space over $k$, we say that $x$ is {\em
$k$-proportional} to $y$ if $x=cy$ for some nonzero scalar $c$.

\subsubsection{Definition}\label{subsubsection:SegreDefinition} 
A {\em Segre matrix} $X$ is an object with
the following properties:
\begin{itemize}
\item $X$ is an $n$ by $n$ matrix with entries in $L$.
\item $X$ is of rank $\leq 1$.
\item $X$ is $k$-general.
\end{itemize} 
If we need to draw attention to the
basic data we say that $X$ is of {\em type
$(k,n,A,L)$}. 

\subsubsection{Key propertries}
Let $X$ be a Segre matrix.
The following  hold:
\begin{itemize}
\item Any matrix with
entries in $A$ to which $X$ is $k$-equivalent is a Segre
matrix.
\item  The transpose $X^T$ is a Segre matrix.
\item There exists a
factorization
$X=uv$ where
$u$ (resp.~$v$) is a column (resp.~row) vector with
entries in the fraction field of $A$. \item Given any
such factorization
$X=uv$, the entries of
$u$ (resp.~$v$) are \linebreak $k$-linearly independent.
\item Given any two such factorizations $X=uv=u'v'$,
there exists  unique nonzero $f$ in the fraction
field of $A$ such that
$u'=fu$ and
$v'=f^{-1}v$. 
\end{itemize}
The proofs of these facts are very easy and therefore omitted.
We take these facts for granted in all subsequent work with Segre
matrices.

\subsubsection{Goal}
We aim to put the $k$-equivalence
classes of Segre matrices into explicit bijective
correspondence with the points of an explicitly defined
projective algebraic variety over
$k$. Our
results are summarized by
Theorem~\ref{Theorem:ProjectiveParameterization} below.

\subsection{Examples of Segre 
matrices involving elliptic functions} 
\label{subsection:EllipticSegre}

\subsubsection{The spaces $\Sigma_N$}\label{subsubsection:EllipticSegre}
For background on elliptic functions, see
\cite[Chap.~XX]{WhittakerWatson}. Let
$\sigma(z)$ be the Weierstrass
$\sigma$-function attached to a lattice
$\Lambda\subset\CC$. By construction $\sigma(z)$ has simple zeroes on
the lattice $\Lambda$ and no other zeroes. For each
nonnegative integer $N$, let
$\Sigma_N$ be the space of entire
functions
$f(z)$ such that the meromorphic function
$f(z)/\sigma(z)^N$ is $\Lambda$-periodic.
We have
$\Sigma_0=\CC$
and for $N>0$ we have
$\dim_\CC\Sigma_N=N$. 

\subsubsection{Specification of a type}
Fix an integer $n\geq 2$. 
Clearly the $\CC$-algebra $\bigoplus_{\ell=0}^\infty\Sigma_{2n\ell}$ is
without zero-divisors. It can be shown  that the $\CC$-algebra
$\bigoplus_{\ell=0}^\infty\Sigma_{2n\ell}$ is generated over $\CC$
by $\Sigma_{2n}$. It follows that the quadruple 
\begin{equation}\label{equation:SegreExampleType}
\left(\CC,n,\bigoplus_{\ell=0}^\infty\Sigma_{2n\ell},\Sigma_{2n}
\right)
\end{equation}
is a type. We are going to classify Segre
matrices of this type.

\subsubsection{An analytic construction of Segre matrices}
Let $\vec{\sigma}(z)$ be a row vector of length $n$
with entries forming a $\CC$-basis of $\Sigma_n$,
e.~g.
$$\vec{\sigma}(z)=
\left[\begin{array}{ccccc}
\sigma(z)^n\;\;&\sigma(z)^n\wp(z)\;\;&\sigma(z)^n
\wp'(z)\;\;&
\cdots\;\;&\sigma(z)^n\wp^{(n-2)}(z)
\end{array}\right],
$$
where $$\wp(z)=-\frac{d^2}{dz^2}\log \sigma(z)$$ is the
Weierstrass
$\wp$-function attached to the lattice
$\Lambda$.
It is not difficult to prove that for each 
$t\in
\CC$
the $n$ by $n$ 
matrix 
\begin{equation}\label{equation:AnalyticFactorization}
\vec{\sigma}(z-t/n)^T\vec{\sigma}(z+t/n)
\end{equation}
of entire functions of $z$ is a Segre matrix of type
(\ref{equation:SegreExampleType}) the
$\CC$-equivalence class of which depends only on
$t\bmod{\Lambda}$, not on the choice of $\vec{\sigma}(z)$.
\begin{Proposition}\label{Proposition:EllipticSegre}
 The map sending $t\in \CC$
to the corresponding Segre matrix 
of the form \textup{(\ref{equation:AnalyticFactorization})}
puts the complex torus
$\CC/\Lambda$ in bijective correspondence with the family
of
$\CC$-equivalence classes of Segre matrices
of type \textup{(\ref{equation:SegreExampleType})}.
\end{Proposition}
\proof
In any given fundamental domain
for
$\Lambda$  a not-identically-vanishing function belonging to the space
$\Sigma_n$ has exactly
$n$ zeroes and moreover these zeroes sum to an element of 
$\Lambda$.  Further, the family of functions
$\Sigma_n$ has no zero in common. Now let
$X(z)$ be a matrix $\CC$-equivalent
to a matrix of the form
(\ref{equation:AnalyticFactorization}). Then the functions
in any given row of $X(z)$  have in any given
fundamental domain for $\Lambda$ exactly
$n$ common zeroes, and these must sum to
$t$ modulo $\Lambda$.
Therefore the correspondence in question is
one-to-one.

Now fix
a Segre matrix
$X(z)$ of type (\ref{equation:SegreExampleType}) arbitrarily. Choose in
$X(z)$ a column
$u(z)$ and a row $v(z)$, each with $\CC$-linearly independent entries.
Let
$f(z)$ be the entry of $X(z)$ common to $u(z)$ and $v(z)$; then $f(z)$
does not vanish identically. Let
$p$ (resp.~$q$) be the number of common zeroes of
 the entries of $u(z)$ (resp.~$v(z)$) in any given
fundamental domain for
$\Lambda$. 
Any
subspace of
$\Sigma_{2n}$ defined by prescribing $n+1$ zeroes in a
given fundamental domain for
$\Lambda$ is
$(n-1)$-dimensional over $\CC$; this is a consequence of Riemann-Roch in
genus one. Since the entries of
$u(z)$ (resp.~$v(z)$) are $\CC$-linearly independent, it follows that
$p\leq n$ (resp.~$q\leq n$). Since every
two by two submatrix of
$X(z)$ has vanishing determinant, we  have
$$X(z)=u(z)v(z)/f(z),$$
hence $f(z)$ divides every entry of the product
$u(z)v(z)$, hence $p+q\geq 2n$, hence
$p=q=n$, and hence 
the common zeroes of the entries of the matrix
$u(z)v(z)$ coincide with the zeroes of $f(z)$. 

For some complex numbers
$P_1,\dots,P_{2n}$ summing to $0$ and some nonzero complex number
$C$ we have a factorization
$$f(z)=C\sigma(z-P_1)\cdots \sigma(z-P_{2n}).$$
Moreover, by re-indexing the points
$P_1,\dots,P_{2n}$ if necessary, we can arrange for the
vectors
$$u^*(z)=\frac{u(z)}{C\sigma(z-P_1)\cdots
\sigma(z-P_n)},\;\;\;
v^*(z)=\frac{v(z)}{\sigma(z-P_{n+1})\cdots
\sigma(z-P_{2n})}$$
to have entries that are entire functions of $z$. Now
put
$$t=-(P_1+\cdots+P_n)=P_{n+1}+\cdots+P_{2n}.$$ 
Then the entries of the
vector 
 $u^*(z+t/n)$ (resp.~$v^*(z-t/n)$) belong to
$\Sigma_n$ and hence, since $\CC$-linearly independent,
form a $\CC$-basis for $\Sigma_n$. It
follows that for some
$n$ by $n$ matrices $\Phi$ and $\Psi$ with entries in
$\CC$ we have
$$u^*(z+t/n)=\Phi\vec{\sigma}(z)^T,\;\;\;
v^*(z-t/n)=\vec{\sigma}(z)\Psi,\;\;\;\det \Phi\neq 0,\;\;\;
\det \Psi\neq 0.$$
Finally, we
have
$$X(z)=\Phi\vec{\sigma}(z-t/n)^T\vec{\sigma}(z+t/n)\Psi,$$ i.~e.,
$X(z)$ is $\CC$-equivalent to a Segre matrix of the form
(\ref{equation:AnalyticFactorization}). Therefore the correspondence in
question is onto. \qed

\subsection{An {\em ad hoc} tensor formalism}
\label{subsection:TensorFormalism}   

\subsubsection{Definition of $A^{\otimes \ZZ}$}
\label{subsubsection:InfiniteTensorPowerDef}

Let $A^{\otimes \ZZ}$ be the tensor product over $k$
of copies of $A$ indexed by $\ZZ$, formed
according to the definition
\cite[pp.\ 301-303]{JacquetLanglands}.
By that definition the $k$-algebra
$A^{\otimes \ZZ}$ is generated by symbols of the form
$$\bigotimes_{i\in \ZZ}a_i\;\;\;
(a_i\in A,\;\;a_i=1\mbox{ for $|i|\gg 0$})$$
subject to obvious relations. 
Put
$$a^{(\ell)}:=\bigotimes_{i\in \ZZ}\left\{\begin{array}{rl}
a&\mbox{if $i=\ell$}\\
1&\mbox{if $i\neq \ell$}
\end{array}\right.
$$
for all $a\in A$ and $\ell\in \ZZ$.
More generally, given a matrix
$X$ with entries in $A$ and $\ell\in \ZZ$ we define a
matrix $X^{(\ell)}$ with entries in $A^{\otimes \ZZ}$ by
the rule
$$\left(X^{(\ell)}\right)_{ij}=
\left(X_{ij}\right)^{(\ell)}.$$
 For any subset
$I\subseteq
\ZZ$, let 
$A^{\otimes I}$
denote the $k$-subalgebra of $A^{\otimes \ZZ}$
generated by all elements of the form $a^{(\ell)}$
where $a\in A$ and $\ell\in I$. 
If $I$ is a finite subset of $\ZZ$,
then the
$k$-algebra
$A^{\otimes I}$ can naturally be identified with the usual tensor product
over $k$ of copies of $A$ indexed by $I$.  For any subset
$I\subset\ZZ$, the
$k$-algebra
$A^{\otimes I}$ is characterized in the category of commutative
$k$-algebras with unit by a universal property we need not belabor. We
make the identifications
$$a^{(0)}=a$$
for all $a\in A$, thus equipping $A^{\otimes
\ZZ}$ with the structure of $A$-algebra.
Finally, note that the
$k$-algebra
$A^{\otimes
\ZZ}$ is without zero divisors; this fact plays an extremely important
role in the sequel.

\subsubsection{Partial specializations of $A^{\otimes
\ZZ}$}  Let
$S$ be the set of $k$-algebra homomorphisms $A\rightarrow k$. For all
$a\in A$ and
$s\in S$ we denote the value in
$k$ of $a$ at $s$ by  $a\vert_s$. More generally, given a
matrix
$X$ with entries in
$A$, we define a matrix
$X\vert_s$ with entries in $k$ by the rule
$$(X_{ij})\vert_s=(X\vert_s)_{ij}.$$ Now let
$I$ be a set of integers and let
$$\sbold=\{s_\ell\}_{\ell\in I}\in S^{I}$$ be any family
of points of $S$ indexed by
$I$. We define the {\em partial specialization}
$$\left(a\mapsto a\Vert_\sbold\right):A^{\otimes \ZZ}\rightarrow
A^{\otimes(\ZZ\setminus I)}$$ 
associated to the family $\sbold$ to be the unique
$k$-algebra homomorphism such that
$$a^{(\ell)}\Vert_\sbold=\left\{\begin{array}{rl}
a\vert_{s_\ell}&\mbox{if $\ell\in I$}\\\\
a^{(\ell)}&\mbox{if $\ell\not\in I$}
\end{array}\right.$$
for all $a\in A$ and $\ell\in \ZZ$. 
More generally, given a matrix
$X$ with entries in $A^{\otimes \ZZ }$, we
define a matrix
$X\Vert_\sbold$ with entries in $A^{\otimes (\ZZ\setminus I)}$
by the rule
$$(X\Vert_\sbold)_{ij}= (X_{ij})\Vert_\sbold.$$

\subsubsection{Derangements and their action on $A^{\otimes \ZZ}$}
A {\em
derangement} is by definition a map of the set
$\ZZ$ of integers to itself.  The {\em support} of a derangement
$\sigma$ is by definition the set
$$\{\ell\in
\ZZ\mid
\sigma \ell\neq \ell\},$$ i.~e., the set of integers 
 actually moved by $\sigma$. 
For
each derangement
$\sigma$ we define
$$\sigma_*:A^{\otimes \ZZ}\rightarrow A^{\otimes \ZZ}$$
to be the unique $k$-algebra homomorphism such that
$$\sigma_*(a^{(\ell)})=a^{(\sigma \ell)}$$
for all $a\in A$ and $\ell\in \ZZ$ and more generally, given
any matrix
$Z$ with entries in $A^{\otimes
\ZZ}$, we define a matrix $\sigma_* Z$ with entries in $A^{\otimes \ZZ}$
by the rule
$$(\sigma_* Z)_{ij} =\sigma_*(Z_{ij}).$$
For any derangements $\sigma$ and $\tau$ we have
$$\sigma_*\tau_*=(\sigma \tau)_*.$$

\subsubsection{The bar operation}
We define the {\em bar operation}
$$\left(a\mapsto \bar{a}\right):A^{\otimes \ZZ}\iso
A^{\otimes \ZZ}$$
to be the unique $k$-algebra automorphism 
such that
$$\overline{a^{(\ell)}}=a^{(-\ell)}$$
for all $a\in A$ and $\ell\in \ZZ$. 
 More generally, given a
matrix
$Z$ with entries in $A^{\otimes \ZZ}$ we define
$\overline{Z}$ by the rule
$$(\,\overline{Z}\,)_{ij}=\overline{Z}_{ij}.$$
The bar operation is none other than the
automorphism of
$A^{\otimes \ZZ}$ associated to the sign-reversing
derangement $\ell\mapsto -\ell$. Since
$$\overline{a}=\overline{a^{(0)}}=a^{(0)}=a$$ for all
$a\in A$, the bar operation is an $A$-algebra
automorphism.

\subsubsection{Special derangements}
Given distinct integers
$i$ and
$j$, let $\langle i\mid j\rangle$ be the unique
derangement with support $\{i,j\}$; in other words,
$\langle i\mid j\rangle$ exchanges $i$ and $j$ and
fixes all other integers. Given integers
$i$ and
$j$ (possibly not distinct), let $[i\mapsto j]$ be
the unique derangement with support contained in the set $\{i\}$ sending
$i$ to $j$; in other words, $[i\mapsto j]$ maps $i$ to $j$
and fixes all other integers.

\subsection{Criteria for 
$k$-generality}  
\begin{Lemma}\label{Lemma:GeneralType}
Let integers $\ell_1<\dots<\ell_N$ be given.
Let $u$ be a column vector of length $N$ with entries
in $A$ and put 
$$U=\left[\begin{array}{ccc}
u^{(\ell_1)}&\dots&u^{(\ell_N)}\end{array}
\right]$$
thereby defining an $N$ by $N$ matrix with entries in
$A^{\otimes \{\ell_1,\dots,\ell_N\}}$. The entries of
the vector
$u$ are
$k$-linearly dependent if and only if the determinant
of the matrix $U$ vanishes identically.
\end{Lemma}
\proof 

($\Rightarrow$) By row operations leaving the
determinant unchanged we can transform $U$
to a  matrix with an identically vanishing row.

$(\Leftarrow)$ We proceed by induction on $N$. The case
$N=1$ is trivial; assume now that $N>1$. Without loss of
generality we may assume that the determinant of every
$N-1$ by
$N-1$ submatrix of $U$ is nonvanishing, for otherwise we
are done by induction on $N$. 
Also without loss of generality we may assume that
$(\ell_1,\dots,\ell_N)=(0,\dots,N-1)$. Expanding
$U$ by minors of the first column we obtain a
relation
$$a_1u_1+\cdots+a_Nu_N=0\;\;\;\left(
a_1,\dots,a_N\in A^{\otimes
\{1,\dots,N-1\}},\;\;a_1\cdots a_N\neq 0\right)$$ among
the entries of
$u$. By the Nullstellensatz
there exists
$$\sbold=\left(s_1,\dots,s_{N-1}\right)\in
S^{\{1,\dots,N-1\}}$$ such that 
$$(a_1\cdots a_N)\Vert_\sbold=(a_1\Vert_\sbold)
\cdots (a_n\Vert_\sbold)\neq
0$$ and for any such $\sbold$ we obtain by partial
specialization a nontrivial
$k$-linear relation
$$
\left(a_1u_1+\cdots+a_N
u_N\right)\Vert_\sbold
=\left(a_1\Vert_\sbold\right)\cdot u_1+\cdots
+\left(a_N\Vert_\sbold\right)\cdot u_N=0\
$$  among the entries of $u$. 
\qed

\begin{Lemma}\label{Lemma:GeneralTypePrime}
Let integers $\ell_1<\dots<\ell_n$
be given. Let $X$ be
an
$n$ by
$n$ matrix with entries in
$A$ of rank $\leq 1$. The matrix $X$ is $k$-general if and only if
$\det\left(\sum_{\nu=1}^n
X^{(\ell_\nu)}\right)\neq 0$.
\end{Lemma}
\proof We may assume without loss of generality
that
$(\ell_1,\dots,\ell_n)=(1,\dots,n)$. Let
$u$ be the $j^{th}$  column of
$X$ and let $v$ be the $i^{th}$ row of $X$, where $i$ and $j$
are presently to be chosen in a useful way. Put
$$U=\left[\begin{array}{ccc}
u^{(1)}&\cdots&u^{(n)}\end{array}\right],\;\;\;
V=\left[\begin{array}{c}
v^{(1)}\\
\vdots\\
v^{(n)}\end{array}\right].$$
We claim that
\begin{equation}\label{equation:SegreTool}
\det U\cdot \det
V=\det\left(X^{(1)}+\cdots+X^{(n)}\right)\cdot
X_{ij}^{(1)}\cdots X_{ij}^{(n)}.
\end{equation}
In any case we have
$$X_{ij}X=uv$$ since $X$ is of rank $\leq 1$.
If $X_{ij}=0$, then both sides of (\ref{equation:SegreTool})
vanish. Suppose now that $X_{ij}\neq 0$. 
After localizing $A$ suitably,
we may assume that
$X_{ij}$ is a unit of $A$. Then we have
$$X^{(1)}+\cdots+X^{(n)}=u\;\diag
\left(X_{ij}^{(1)},\dots,X_{ij}^{(n)}\right)^{-1}v,$$
whence (\ref{equation:SegreTool})
after taking determinants on both sides
and rearranging the resulting identity. The claim is  proved.

To prove the implication $(\Rightarrow)$, we
choose $i$ and $j$ so that the entries of
$u$ are
$k$-linearly independent and the entries of $v$ are
$k$-linearly independent. Then the entry
$X_{ij}$ common to $u$ and $v$ is nonzero and the
left side of (\ref{equation:SegreTool}) is nonvanishing
by Lemma~\ref{Lemma:GeneralType}. It follows that
the determinant in question does not vanish. 

To prove the implication $(\Leftarrow)$,  we choose $i$
and $j$ so that 
$X_{ij}\neq 0$. Then the right side of (\ref{equation:SegreTool}) 
is nonvanishing and hence
neither
$\det U$ or
$\det V$ vanish. By Lemma~\ref{Lemma:GeneralType} the
entries of
$u$ are $k$-linearly independent and the entries of $v$
are $k$-linearly independent. It follows that the matrix $X$ is
$k$-general.
\qed

\begin{Proposition}\label{Proposition:GeneralTypeTilde}
Let $X$ be an $n$ by $n$ matrix with entries in $A$
of rank $\leq 1$. The following conditions are
equivalent:
\begin{equation}\label{equation:GT1}
\textup{\mbox{$X$ is $k$-general.}}
\end{equation}
\begin{equation}\label{equation:GT3}
\Delta_{\ell=1}^{n+1}X^{(\ell)}\neq 0.
\end{equation}
\begin{equation}\label{equation:GT6}
\abel_{\ell=0}^{n+1}X^{(\ell)}\neq 0.
\end{equation}
\end{Proposition}
\proof For distinct
$i,j\in
\{0,\dots,n+1\}$ put
$$D_{ij}=\det\left(\sum_{\ell\in
\{0,\dots,n+1\}\setminus\{i,j\}}X^{(\ell)}\right).
$$
Consider two more conditions:
\begin{equation}\label{equation:GT7}
D_{ij}\neq 0\;\;\textup{\mbox{for some distinct $i,j\in
\{0,\dots,n+1\}$.}}
\end{equation}
\begin{equation}\label{equation:GT8}
D_{ij}\neq 0\;\;\textup{\mbox{for all distinct $i,j\in
\{0,\dots,n+1\}$.}}
\end{equation}
We
have implications
$$(\ref{equation:GT8})\Rightarrow(\ref{equation:GT3})
\Rightarrow(\ref{equation:GT7})
\Rightarrow
(\ref{equation:GT1})\Rightarrow(\ref{equation:GT8})
\Rightarrow(\ref{equation:GT6}),
$$
the first two by definition of the discriminant,
the next two by Lemma~\ref{Lemma:GeneralTypePrime},
and the last by identity
(\ref{equation:RankOneDiagonalAbeliantEvaluation}). To complete
the proof it suffices to prove the implication 
(\ref{equation:GT6})$\Rightarrow$(\ref{equation:GT7}). 
Put 
$Z=\abel_{\ell=0}^{n+1}X^{(\ell)}$ to abbreviate notation.
By hypothesis there exist indices $i$ and $j$ such that 
$Z_{ij}\neq 0$.
We then have
$$
0\neq Z_{ij}\cdot\langle 0\mid n+1\rangle_* 
Z_{ij}=
Z_{ij}Z_{ji}=Z_{ii}Z_{jj}
=
D_{0i}D_{i,n+1}D_{0j}D_{j,n+1},$$
the first equality by 
identity
(\ref{equation:AbeliantTransposeSymmetry}),
the second by fact (\ref{equation:RankOneInRankOneOut})
and the third by identity
(\ref{equation:RankOneDiagonalAbeliantEvaluation}).
Therefore condition (\ref{equation:GT6})
does indeed imply condition (\ref{equation:GT7}).
\qed

\subsection{Normalization and self-similarity}
\subsubsection{Definitions}
Let an $n$ by $n$ matrix
$X$ with entries in $A$ of rank $\leq 1$ be given.
Let
$$\sbold=(s_1,\dots,s_{n+1})\in
S^{\{1,\dots,n+1\}}$$ 
be given.
We say that $X$ is {\em
$\sbold$-normalized} under the following conditions:
\begin{equation}\label{equation:NormalizationTwo}
\mbox{For $i,j,\ell=1,\dots,n$ we have $X_{ij}\vert_{s_\ell}\neq 0$ if
and only if $i=j=\ell$.}
\end{equation}
\begin{equation}\label{equation:NormalizationOne}
\mbox{For $i,j=1,\dots,n$ we have $X_{ij}\vert_{s_{n+1}}=1$.}
\end{equation}
 We say that
$X$ is {\em
$\sbold$-self-similar} under the following conditions:
\begin{equation}
\label{equation:SelfSimilarOne}
\Delta\left(X\vert_{s_1},\dots,X\vert_{s_{n+1}}\right)\neq
0.
\end{equation}
\begin{equation}\label{equation:SelfSimilarTwo}
\mbox{$\abel
\left(X,X\vert_{s_1},\dots,X\vert_{s_{n+1}}\right)$ is
$k$-proportional to $X$.}
\end{equation}
\begin{Proposition}
\label{Proposition:GeneralTypeAfterThought}
Fix an $n$ by $n$ matrix
$X$ with entries in $A$ of rank $\leq 1$. Fix
$\sbold=(s_1,\dots,s_{n+1})\in S^{\{1,\dots,n+1\}}$.  (i) If
$X$ is
$\sbold$-normalized then condition
\textup{(\ref{equation:SelfSimilarOne})} holds. (ii) If condition
\textup{(\ref{equation:SelfSimilarOne})} holds
then $X$ is $k$-general.
\end{Proposition}
\proof (i) This follows from identity
(\ref{equation:DiscriminantSpecialCase}).
(ii) Put $\Delta=\Delta_{\ell=1}^{n+1}X^{(\ell)}$. By hypothesis
$\Delta\Vert_\sbold\neq 0$, hence
$\Delta\neq 0$, whence the result by
Proposition~\ref{Proposition:GeneralTypeTilde}.
\qed
\begin{Proposition}\label{Proposition:NormalizationOne}
Fix a Segre matrix $X$ and
$\sbold=(s_1,\dots,s_{n+1})\in S^{\{1,\dots,n+1\}}$.
There exist $n$ by $n$ matrices $\Phi$ and
$\Psi$ with entries in $k$ such
that
$$\Phi X\Psi=\abel\left(X,X\vert_{s_1},\dots,
X\vert_{s_{n+1}}\right),\;\;\;
(\det
\Phi)^2(\det
\Psi)^2=\Delta
\left(X\vert_{s_1},\dots,X\vert_{s_{n+1}}\right).
$$ 
Moreover if $X$ is $\sbold$-normalized then $\Phi$ and $\Psi$
may be taken diagonal.
\end{Proposition}
\proof This follows directly from
key relation (\ref{equation:KeyRelation}). 
\qed
\begin{Lemma}\label{Lemma:NormalizationTwo}
Fix a Segre matrix $X$ and
$\sbold=(s_1,\dots,s_{n+1})\in S^{\{1,\dots,n+1\}}$.
If $X$ is $\sbold$-self-similar then there exist diagonal
matrices $\Phi$ and $\Psi$
with entries in $k$ such
that $\det
\Phi \cdot \det\Psi\neq 0$ and $\Phi^{-1} X\Psi^{-1}$ is
$\sbold$-normalized. 
\end{Lemma}
\proof 
By hypothesis (\ref{equation:SelfSimilarTwo}) there exists a nonzero
scalar
$c$ such that
$$X=c\cdot\abel\left(X,X\vert_{s_1},
\dots,X\vert_{s_{n+1}}\right).$$
It follows by
identity (\ref{equation:AbeliantDuplicateVanishing}) that for
$\ell=1,\dots,n$ every entry of  the matrix
$$X\vert_{s_\ell}=c\cdot\abel\left(X\vert_{s_\ell},
X\vert_{s_1},\dots,X\vert_{s_{n+1}}\right)
$$
vanishes save possibly the $\ell^{th}$ diagonal entry.
Let $E$ be the $n$ by $n$ matrix with all entries equal to $1$.
Since $X\vert_{s_{n+1}}$ is of rank $\leq 1$ we can write 
$$X\vert_{s_{n+1}}=\Phi E \Psi$$
where $\Phi$ and $\Psi$ are diagonal matrices with entries in $k$. 
By hypothesis (\ref{equation:SelfSimilarOne}) and
identity (\ref{equation:DiscriminantSpecialCase})
neither is it possible for $\det \Phi\cdot \det\Psi$ to vanish,  nor for
there to exist some index
$\ell=1,\dots,n$ such that the $\ell^{th}$ diagonal entry of
$X\vert_{s_\ell}$ vanishes. Therefore the pair
$(\Phi,\Psi)$ has all the desired properties. 
\qed

\begin{Proposition}\label{Proposition:AffineCovering}
Fix a Segre matrix $X$ and $\sbold\in
S^{\{1,\dots,n+1\}}$. Put
$$Z=\abel_{\ell=0}^{n+1}X^{(\ell)},\;\;\;
\Delta=\Delta_{\ell=1}^{n+1}X^{(\ell)}.$$
Assume that $\Delta\Vert_\sbold\neq 0$.
(i) Up to $k$-proportionality $Z\Vert_\sbold$  is the unique
$\sbold$-self-similar Segre matrix
$k$-equivalent to $X$. (ii) There exists
a unique 
$\sbold$-normalized Segre matrix
$k$-equivalent to
$X$.
\end{Proposition}
\proof (i) By
Proposition~\ref{Proposition:NormalizationOne} the matrix
$$Z\Vert_\sbold=
\abel\left(X,X\vert_{s_1},\dots,X\vert_{s_{n+1}}\right)
$$
is a Segre matrix
$k$-equivalent to $X$. Further, $Z\Vert_\sbold$ is
$\sbold$-self-similar by the iterated abeliant
identity (\ref{equation:IteratedAbeliants}). Finally, any two
$k$-equivalent
$\sbold$-self-similar Segre matrices are
$k$-proportional by the abeliant transformation law
(\ref{equation:FundamentalLaw}). 

(ii) We may assume without
loss of generality that
$X$ is
$\sbold$-self-similar. By Lemma~\ref{Lemma:NormalizationTwo} 
there exists at least one
$\sbold$-normalized Segre matrix $k$-equivalent to
$X$.
Now suppose that
$Y$ and $Y'$ are
$\sbold$-normalized Segre matrices both \linebreak
$k$-equivalent to
$X$. By
Proposition~\ref{Proposition:NormalizationOne} and our additional
assumption that
$X$ is
$\sbold$-self-similar, we have
$X=\Phi Y\Psi=\Phi'Y'\Psi'$ where
$\Phi$, $\Psi$, $\Phi'$, $\Psi'$ are nonsingular diagonal
matrices with entries in $k$, hence we have 
$\Phi E \Psi=\Phi'E\Psi'$
where $E$ is the $n$ by $n$ matrix with all entries equal to
$1$, hence there exists a nonzero scalar $c$ such that
$\Phi'=c\Phi$ and $\Psi'=c^{-1}\Psi$, and hence $Y=Y'$.
\qed

\subsection{The abstract Abel map: definition and
key properties}

\subsubsection{Definition}
The {\em abstract Abel map}
by definition sends each Segre matrix $X$ to the $n$ by $n$ matrix
$\abel_{\ell=0}^{n+1}X^{(\ell)}$ with entries in $A^{\otimes
\{0,\dots,n+1\}}$.  The abstract Abel map generalizes and abstracts the
explicit algebraic representation of the Abel map studied
in \cite{Anderson}. 

\subsubsection{Catalog of key
properties}\label{subsubsection:BasicProperties}  Fix a Segre
matrix
$X$.  Let 
$$Z=\abel_{\ell=0}^{n+1}X^{(\ell)},\;\;\;
\Delta=\Delta_{\ell=1}^{n+1}X^{(\ell)}$$
be the image of $X$ under the
abstract Abel map and the naturally associated discriminant,
respectively.  The following hold:
\begin{itemize}
\item $Z\neq 0$ and $\Delta\neq 0$, by
Proposition~\ref{Proposition:GeneralTypeTilde}.\\
\item If $X$ is replaced by a $k$-equivalent matrix,
then $Z$ and $\Delta$ are replaced by nonzero scalar
multiples, by (\ref{equation:FundamentalLaw}) 
and (\ref{equation:DiscriminantTransformationLaw}).\\
\item If $X$ is replaced by $X^T$, then
$Z$ is replaced
by $Z^T$ and $\Delta$ remains unchanged, by
(\ref{equation:AbeliantTransposeSymmetry}) and
(\ref{equation:DiscriminantTransposeSymmetry}).\\
\item $\langle 0\mid n+1\rangle_*Z=Z^T$, by
(\ref{equation:AbeliantTransposeSymmetry}).\\
\item $\pi_*
Z_{ij}=Z_{\pi i, \pi j}$
for any bijective
derangement $\pi$ supported in
$\{1,\dots,n\}$, \linebreak by (\ref{equation:AbeliantSymmetry}).\\
\item $[1\mapsto 2]_*Z_{12}=Z_{11}$, by
(\ref{equation:AbeliantDegeneration}).\\
\item $Z$ is of rank $\leq 1$, by (\ref{equation:RankOneInRankOneOut}).\\
\item $\displaystyle\Delta=
\langle
2|n+1\rangle_*\langle 0|1\rangle_* 
 Z_{11}\cdot
\langle 0|1\rangle_* Z_{11}\cdot
\langle 0|2\rangle_* Z_{22}\cdot
\prod_{\ell=3}^n
\left(\langle 0|\ell\rangle_*
Z_{\ell\ell}\right)^2$, \linebreak
by (\ref{equation:DiscriminantDefinitionReformulation}).\\
\item $\abel_{\ell=0}^{n+1}[0\mapsto
-\ell]_*Z=\Delta \cdot \overline{Z}$, by
(\ref{equation:IteratedAbeliants}). \\
\item $\displaystyle Z_{ij}\in \mbox{$k$-span of }
L\;\cdot\;\;
\prod_{b\in \{1,\dots,n\}\setminus\{j\}}L^{(b)}
\;\;\cdot \;L^{(n+1)}\;\cdot\;\;\prod_{a\in
\{1,\dots,n\}\setminus \{i\}}L^{(a)}
$, \linebreak
by (\ref{equation:FourPermutationExpansion}).
\end{itemize}

\begin{Proposition}[``The abstract
Abel theorem'']
\label{Proposition:AbelSeparation} Let
$X$ and $X'$ be \linebreak Segre matrices with corresponding images $Z$
and
$Z'$ under the abstract Abel map, respectively. 
Then $X'$ is $k$-equivalent to $X$ 
if and only if
$Z'$ is \linebreak $k$-proportional to $Z$.

\end{Proposition}
\proof 

($\Rightarrow$) This follows directly from identity
(\ref{equation:FundamentalLaw}).

($\Leftarrow$) Put
$$\Delta=\Delta_{\ell=1}^{n+1}X^{(\ell)},\;\;\;
\Delta'=\Delta_{\ell=1}^{n+1}(X')^{(\ell)}.$$
By Proposition~\ref{Proposition:GeneralTypeTilde} neither $\Delta$
nor $\Delta'$ vanish identically. By the Nullstellensatz there
exists
$$\sbold=(s_1,\dots,s_{n+1})\in S^{\{1,\dots,n+1\}}$$
such that 
$$\Delta\Vert_\sbold\neq 0,\;\;\;\Delta'\Vert_\sbold\neq 0.$$
By 
Proposition~\ref{Proposition:NormalizationOne} and hypothesis
we have
$$X\sim Z\Vert_\sbold\sim Z'\Vert_\sbold\sim
X',$$
where $\sim$ denotes $k$-equivalence.
\qed

\subsection{Characterization of the image of the abstract
Abel map}

\subsubsection{$J$-matrices}
A {\em $J$-matrix} $Z$ 
is by definition an object with the following
properties:
\begin{equation}\label{equation:JMat0}
\mbox{$Z$ is an $n$ by $n$ matrix with entries in
the $k$-span of $L\cdot A^{\otimes \{1,\dots,n+1\}}$.}
\end{equation}
\begin{equation}\label{equation:JMat1}
Z\neq 0.
\end{equation}
\begin{equation}\label{equation:JMat2}
\mbox{$Z$ is of rank $\leq 1$.}
\end{equation}
\begin{equation}\label{equation:JMat3}
\abel_{\ell=0}^{n+1}[0\mapsto -\ell]_*Z=\Delta \cdot
\overline{Z}\;\;\mbox{for some $0\neq \Delta\in A^{\otimes
\{1,\dots,n+1\}}$}.
\end{equation}
Since $A^{\otimes \ZZ}$ is a $k$-algebra
without zero-divisors,
$Z$ uniquely determines
$\Delta$.  We call
$\Delta$ the {\em discriminant} of the $J$-matrix $Z$.
If we need to call attention to the basic data we say that
$Z$ is a $J$-matrix of {\em type $(k,n,A,L)$}. In view of the
properties catalogued in
\S\ref{subsubsection:BasicProperties}, it is clear that
every matrix in the image of the abstract Abel map is
automatically a
$J$-matrix.

\begin{Proposition}[``The abstract Jacobi
inversion theorem'']\label{Proposition:Clincher} 
Fix a \linebreak $J$-matrix $Z$ with associated discriminant
$\Delta$. (i) For all 
$\sbold
\in S^{\{1,\dots,n+1\}}$ the partial specialization
$\Delta\Vert_\sbold$ is a scalar and for some $\sbold$
that scalar does not vanish. (ii) For any
$\sbold
\in S^{\{1,\dots,n+1\}}$ such that 
the scalar $\Delta\Vert_\sbold$
does not vanish, the corresponding partial specialization
$Z\Vert_\sbold$ is a Segre matrix with image
under the abstract Abel map $k$-proportional to $Z$. 
\end{Proposition}
\proof 

(i) This follows from the Nullstellensatz. 

(ii) Put
$$X=Z\Vert_\sbold.$$ The matrix 
$X$ is an
$n$ by $n$ matrix with entries in
$L$ by condition (\ref{equation:JMat0}) and of rank $\leq 1$
by condition (\ref{equation:JMat2}).
By applying firstly the partial specialization 
operation $\Vert_\sbold$ and secondly the bar operation to both sides
of the identity figuring in condition
(\ref{equation:JMat3}), we obtain the relation
$$\abel_{\ell=0}^{n+1}X^{(\ell)}=
\Delta\Vert_\sbold\cdot Z.
$$
The right side  does not
vanish by condition (\ref{equation:JMat1}) combined with our hypothesis
that $\Delta\Vert_{\sbold}\neq 0$. It
follows by Proposition~\ref{Proposition:GeneralTypeTilde} that
$X$ is $k$-general and hence a Segre matrix.
It follows as well that the image of $X$ under the abstract
Abel map is $k$-proportional to $Z$. 
\qed

\subsubsection{Jacobi matrices}
\label{subsubsection:JacobiDefinition}
A {\em Jacobi matrix} $Z$ is by definition an object with
the following properties:
\begin{equation}\label{equation:JacMat0}
\mbox{$Z$ is an $n$ by
$n$ matrix with entries in
$A^{\otimes
\{0,\dots,n+1\}}$.}
\end{equation}
\begin{equation}\label{equation:JacMat1}
Z\neq 0.
\end{equation}
\begin{equation}\label{equation:JacMat2}
Z_{12}\in \mbox{$k$-span of }
L\cdot L^{(1)}\cdot L^{(2)}\cdot
\left(L^{(3)}\right)^2\cdots \left(L^{(n)}\right)^2\cdot
L^{(n+1)}.
\end{equation}
\begin{equation}\label{equation:JacMat3}
Z_{11}=[1\mapsto 2]_*Z_{12}.
\end{equation}
\begin{equation}\label{equation:JacMat4}
\begin{array}{l}
\mbox{$\pi_*
Z_{ij}=Z_{\pi i, \pi j}$
for any bijective derangement $\pi$ supported}\\
\mbox{in
$\{1,\dots,n\}$.}
\end{array}
\end{equation}
\begin{equation}\label{equation:JacMat5}
\left|\begin{array}{cc}
Z_{11}&Z_{12}\\
Z_{21}&Z_{22}\end{array}\right|=0.
\end{equation}
\begin{equation}\label{equation:JacMat6}
\abel_{\ell=0}^{n+1}[0\mapsto
-\ell]_*Z=\Delta\cdot \overline{Z},
\end{equation}
where
\begin{equation}\label{equation:JacMat7}
\Delta=\langle
2|n+1\rangle_*\langle 0|1\rangle_*  Z_{11}\cdot\langle 0|1\rangle_*
Z_{11}\cdot
\langle 0|2\rangle_* Z_{22}\cdot
\prod_{\ell=3}^n
\left(\langle 0|\ell\rangle_* Z_{\ell\ell}\right)^2.
\end{equation}
We call $\Delta$ the {\em discriminant}
of the Jacobi matrix $Z$.  If
we need to draw attention to the basic data we say that $Z$
is of {\em type
$(k,n,A,L)$}. It is clear that the set of $k$-proportionality
classes of Jacobi matrices forms a projective algebraic variety.
Moreover, in view of the properties cataloged in
\S\ref{subsubsection:BasicProperties}, it is clear that the abstract
Abel map  takes its values in the set of Jacobi matrices.
Note that a Jacobi matrix $Z$ is uniquely determined by its entry
$Z_{12}$.

\begin{Lemma}\label{Lemma:Lynchpin}
The discriminant of a Jacobi matrix does not vanish identically
and belongs to
$A^{\otimes\{1,\dots,n+1\}}$. 
\end{Lemma}
\proof Fix a Jacobi matrix $Z$ with associated discriminant $\Delta$.
 By conditions
(\ref{equation:JacMat1}) and (\ref{equation:JacMat4})
either every diagonal entry
of
$Z$ is nonvanishing or every off-diagonal entry of
$Z$ is nonvanishing; but then by condition (\ref{equation:JacMat5})
every entry $Z$ is nonvanishing. By definition
(\ref{equation:JacMat7}) it follows that
$\Delta$ does not vanish identically.
By conditions (\ref{equation:JacMat2}) and (\ref{equation:JacMat3}) we
have
$$Z_{11}\in A^{\otimes
\{0,2\dots,n+1\}}.$$ By condition (\ref{equation:JacMat4}) we have
$$Z_{\ell
\ell}=\langle 1\mid \ell\rangle_* Z_{11}\in A^{\otimes
\{0,\dots,n+1\}\setminus
\{\ell\}}$$ 
and further
$$\langle \ell\mid 0\rangle_*Z_{\ell\ell}\in A^{\otimes
\{1,\dots,n+1\}},\;\;\langle 2\mid n+1\rangle_*\langle 1\mid
0\rangle_* Z_{11}\in A^{\otimes
\{1,\dots,n+1\}}.$$
By definition (\ref{equation:JacMat7}) it
follows that $\Delta$ belongs to
$A^{\otimes
\{1,\dots,n+1\}}$.
\qed

\begin{Proposition}\label{Proposition:JJacEquiv}
Every Jacobi matrix is a $J$-matrix and vice versa.
\end{Proposition}
\proof By Proposition~\ref{Proposition:Clincher} every
$J$-matrix is in the image of the abstract Abel map and
hence a Jacobi matrix in view of the properties cataloged in
\S\ref{subsubsection:BasicProperties}. Thus the ``vice
versa'' part of the proposition is proved. Now fix a Jacobi matrix $Z$
with associated discriminant
$\Delta$. We verify that $Z$ has the properties required of a $J$-matrix
as follows. To abbreviate notation
temporarily let
$V$ denote the $k$-span of $L\cdot A^{\otimes
\{1,\dots,n+1\}}$. We have $Z_{12}\in V$ by condition
(\ref{equation:JacMat2}). We have  $[1\mapsto
2]_*V\subseteq V$ and hence $Z_{11}\in V$ by condition
(\ref{equation:JacMat3}). We have  $\pi_*V\subseteq V$
for any bijective derangement
$\pi$ supported in
$\{1,\dots,n\}$ and hence $Z_{ij}\in V$ for all $i,j\in
\{1,\dots,n\}$ by condition (\ref{equation:JacMat4}). Therefore $Z$
 satisfies condition (\ref{equation:JMat0}).
Conditions (\ref{equation:JMat1})
and (\ref{equation:JacMat1}) are exactly the same. Clearly $Z$
satisfies condition (\ref{equation:JMat2})
by conditions (\ref{equation:JacMat4}) and
(\ref{equation:JacMat5}). 
Finally, $Z$ satisfies condition 
(\ref{equation:JMat3}) by condition (\ref{equation:JacMat6})
combined with Lemma~\ref{Lemma:Lynchpin}.
\qed
\begin{Theorem}
\label{Theorem:ProjectiveParameterization}
(i) The set of
$k$-proportionality classes of Jacobi matrices forms a
projective algebraic variety. (ii) The
abstract Abel map takes values in the set of Jacobi
matrices.  (iii) The abstract Abel map puts the \linebreak
$k$-equivalence classes of Segre matrices
into bijective correspondence with the $k$-proportionality
classes of Jacobi matrices. 
\end{Theorem}
\proof

(i,ii) These facts have already been noted above.

(iii) The correspondence in question is well-defined and one-to-one by
Proposition~\ref{Proposition:AbelSeparation}.
The correspondence is onto
by Propositions~\ref{Proposition:Clincher}
and \ref{Proposition:JJacEquiv}.  \qed

\subsubsection{Remark}
\label{subsubsection:GeometricalPicture}
We briefly describe the big picture in more geometrical language.
We temporarily introduce the following notation:
\begin{itemize}
\item Let $V$ be the quasi-affine variety
of Segre matrices.
\item Let $G$ denote the product of two
copies of the $n$ by $n$ general linear group over $k$
and let $G$ act in the obvious way on $V$.
\item Let $J$ be the projective variety of
$k$-proportionality classes of Jacobi matrices. 
\end{itemize}
 For
each
$\sbold\in S^{\{1,\dots,n+1\}}$:
\begin{itemize}
\item Let $V_\sbold$ be the
open subvariety of $V$ consisting of Segre matrices $X$
satisfying the inequality
$\Delta_{\ell=1}^{n+1}\left(X\vert_{s_\ell}\right)\neq
0$.
\item Let $U_\sbold$ be the affine variety consisting
of
$n$ by $n$ matrices $X$ with entries in $L$ that are
of rank $\leq 1$ and $\sbold$-normalized.
\item Let
$U'_\sbold$ be the quasi-projective variety
consisting of
$\sbold$-self-similar Segre matrices modulo
$k$-proportionality.
\item  Let $J_\sbold$ be the open subvariety of $J$
consisting of points represented by
Jacobi matrices with discriminant $\Delta$ such that
$\Delta\Vert_\sbold\neq 0$.  
\end{itemize}
By
Proposition~\ref{Proposition:GeneralTypeTilde} 
and the Nullstellensatz the
open subvarieties $V_\sbold$ cover $V$. By identity
(\ref{equation:DiscriminantTransformationLaw}) the
variety
$V_\sbold$ is $G$-stable.
By
Proposition~\ref{Proposition:GeneralTypeAfterThought}
the variety $U_\sbold$ is contained in $V_\sbold$.
By Proposition~\ref{Proposition:AffineCovering}  the
quotient $V_\sbold/G$ can naturally be identified with
$U'_\sbold$ and also with $U_\sbold$.  By
Lemma~\ref{Lemma:Lynchpin} and the
Nullstellensatz the open sets
$J_\sbold$ cover $J$. By
Propositions~\ref{Proposition:AffineCovering},
\ref{Proposition:AbelSeparation},
\ref{Proposition:Clincher} and \ref{Proposition:JJacEquiv}, the
set
$J_\sbold$ can naturally be identified with
$U'_\sbold$ and hence also with $U_\sbold$.  
The upshot is that the family $\{U_\sbold\}$ can be viewed as an affine
open covering of $J$.

\subsection{Examples of Jacobi matrices involving 
elliptic functions}\label{subsection:EllipticJacobi}
\subsubsection{The set up}
We continue in the set up of \S\ref{subsubsection:EllipticSegre}.
By  Proposition~\ref{Proposition:EllipticSegre}
combined with
Theorem~\ref{Theorem:ProjectiveParameterization} every
Jacobi matrix of type (\ref{equation:SegreExampleType}) is 
$\CC$-proportional to the image of some matrix of the form
(\ref{equation:AnalyticFactorization}) under the abstract Abel
map for a value of the parameter $t\in \CC$ uniquely
determined modulo the period lattice $\Lambda$. To figure
out what these Jacobi matrices actually look like, we
are going to apply the abstract Abel map to the matrix
(\ref{equation:AnalyticFactorization}). But before we
proceed we have a notational collision to deal
with: $f^{(\ell)}$ denotes the
$\ell^{th}$ derivative of $f$ in the present context.
To fix the problem we think of and write out the
``superscript
$\ell$'' operation defined in
\S\ref{subsubsection:InfiniteTensorPowerDef} as
the high school algebra operation of ``substitution of
the variable $z_\ell$ for the variable
$z$.'' 

\subsubsection{A classical determinant identity and related
abeliant identity} By combining the classical identity
\begin{equation}\label{equation:FrobStick}
\left|\begin{array}{ccccccccccc}
1&\frac{\wp(z_1)}{1!}&-\frac{\wp'(z_1)}{2!}&\cdots
&\frac{(-1)^{n-2}\wp^{(n-2)}(z_1)}{(n-1)!}
\\
\vdots&\vdots&\vdots&&\vdots\\
1&\frac{\wp(z_n)}{1!}&-\frac{\wp'(z_n)}{2!}&
\cdots&
\frac{(-1)^{n-2}\wp^{(n-2)}(z_n)}{(n-1)!}
\end{array}\right|=\frac{\displaystyle
\sigma\left(\sum_{i=1}^{n} z_i\right)
\cdot \prod_{1\leq i<j\leq
n}\sigma(z_i-z_j)} {\displaystyle
\prod_{i=1}^{n}\sigma(z_i)^n}
\end{equation}
(see \cite[p.~179]{FrobStick} or \cite[Chap.\ XX,
Misc.\ Ex.\ 21]{WhittakerWatson}) with abeliant identity
(\ref{equation:KeyRelationCompanion}), we find after a
straightforward calculation that\\
\begin{equation}\label{equation:BigDealAbeliantIdentity}
\begin{array}{rlll}
&&\left(\abel_{\ell=0}^{n+1}
\vec{\sigma}(z_\ell-t/n)^T\vec{\sigma}(z_\ell+t/n)\right)_{ij}\bigg
/\left(\det\left[\begin{array}{c}
\frac{\vec{\sigma}^{(n)}(0)}{n!}\\
\frac{\vec{\sigma}^{(n-2)}(0)}{(n-2)!}\\
\vdots\\
\frac{\vec{\sigma}(0)}{0!}
\end{array}\right]\right)^4\\\\
=&&\displaystyle\sigma
\left(t+\sum_{\ell\in
\{0,\dots,n+1\}\setminus\{0,i\}}z_\ell\right)\times\sigma
\left(t-\sum_{\ell\in
\{0,\dots,n+1\}\setminus\{i,n+1\}}z_\ell\right)\\\\
&\times&\displaystyle\sigma
\left(t+\sum_{\ell\in
\{0,\dots,n+1\}\setminus\{j,n+1\}}z_\ell\right)\times\sigma
\left(t-\sum_{\ell\in
\{0,\dots,n+1\}\setminus\{0,j\}}z_\ell\right)\\\\
&\times &\displaystyle
\prod_{ \begin{array}{c}
\scriptstyle\alpha,\beta\in \{0,\dots,n+1\}\setminus \{0,i\}\\
\scriptstyle\alpha<\beta
\end{array}}\sigma(z_\alpha-z_\beta)
\times
\prod_{\begin{array}{c}
\scriptstyle\alpha,\beta\in \{0,\dots,n+1\}\setminus \{i,n+1\}\\
\scriptstyle\alpha<\beta
\end{array}}\sigma(z_\alpha-z_\beta)
\\\\
&\times&\displaystyle\prod_{\begin{array}{c}
\scriptstyle\alpha,\beta\in \{0,\dots,n+1\}\setminus \{j,n+1\}\\
\scriptstyle\alpha<\beta
\end{array}}\sigma(z_\alpha-z_\beta)
\times\prod_{\begin{array}{c}
\scriptstyle\alpha,\beta\in \{0,\dots,n+1\}\setminus \{0,j\}\\
\scriptstyle\alpha<\beta
\end{array}}\sigma(z_\alpha-z_\beta).
\end{array}
\end{equation}
Identity (\ref{equation:BigDealAbeliantIdentity}) granted,
it is not difficult to verify that the variety of \linebreak
$\CC$-proportionality classes of Jacobi matrices
is a complex manifold isomorphic to the complex torus
$\CC/\Lambda$. We omit further details.

\section{Elementary construction
of Jacobians}\label{section:MatrixRepresentation}  

\subsection{Basic notation and terminology}
As above, let $k$ be an algebraically closed field.
We work in the category of quasi-projective varieties over $k$. As
above, let
$\curve$ be a nonsingular projective algebraic curve of genus $g$.
We assume that $g>0$.
{\em Divisors} are divisors of
$\curve$. Given a divisor $D$, we write 
$L(D)=H^0(\curve,\OO_\curve(D))$ and $\ell(D)=\dim_k L(D)$.

\subsection{Matrix representation of divisor
classes}
\label{subsection:Dictionary}

\begin{Lemma}\label{Lemma:FundamentalUniquenessPrinciple}
Let $D$ be a divisor such that $\deg D\geq 2g$.
Then we have 
$$\ell(D)=\deg D-g+1.$$
Moreover, we have
$$\ell(D')\geq \ell(D)\Rightarrow\deg D'\geq \deg D$$
for all divisors $D'$.
\end{Lemma}
\proof Riemann-Roch.
\qed

\subsubsection{$G$-forms}\label{subsubsection:GForms}
Let $G$ be a divisor such that
$$\deg G\equiv 0\bmod{2},\;\;\;\;\frac{1}{2}\deg G\geq 2g$$ and put 
$$n=\frac{1}{2}\deg
G-g+1.$$
 By
definition a {\em
$G$-form} $X$ is an object with the following
properties:
\begin{itemize}
\item $X$ is an $n$ by $n$ matrix with entries in
$L(G)$.
\item Every two by two submatrix of $X$ has vanishing determinant.
\item There exists in $X$ some row  and also some column with
$k$-linearly independent entries.
\end{itemize}
A $G$-form is the same thing as a
Segre matrix of type
$$\left(k,n,
\bigoplus_{m=0}^\infty L(mG),L(G)\right).$$ If
$G> 0$,  then a $G$-form can also be viewed as a Segre matrix of type
$$\left(k,n,H^0(\curve\setminus \supp
G,\OO_\curve),L(G)\right).$$
Here and below $\supp D$ denotes the support of a divisor $D$.
Our immediate goal is to put the
$k$-equivalence classes of
$G$-forms in canonical bijective correspondence with the
divisor classes of degree $\frac{1}{2}\deg
G$. The ``dictionary'' we ultimately obtain is summarized by
Proposition~\ref{Proposition:Dictionary} below.

\subsubsection{Representation of divisors of degree
$\frac{1}{2}\deg G$ by $G$-forms}
\label{subsubsection:DivMatrixRep}  Let $G$ and $n$ be as above. Let
$D$ be a divisor of degree
$\frac{1}{2}\deg G$. Let
$u$ (resp.~$v$) be a column (resp.~row) vector with
entries forming a
$k$-basis for
$L(D)$ (resp.~$L(G-D)$). Put 
$X=uv$. It is clear that that $X$
is a $G$-form the $k$-equivalence class of which depends only on $D$,
not on the choice of vectors $u$ and $v$.
In this situation we say that the $G$-form
$X$ {\em represents} the divisor
$D$ of degree $\frac{1}{2}\deg G$.
We claim that for
every divisor $D'$ in the divisor class of $D$ and every
$G$-form
$X'$ representing
$D'$, the $G$-form $X$ is $k$-equivalent to $X'$. To
prove the claim, write 
$D=D'+(f)$ where $f$ is a nonzero meromorphic function on
$\curve$. Then 
$$L(D')=f\cdot L(D),\;\;\;L(G-D')=f^{-1}\cdot L(G-D),$$
hence the $G$-form
$X=uv=(fu)(f^{-1}v)$ represents not only $D$ but also $D'$,
and
hence $X$ is $k$-equivalent to
$X'$. The claim is proved.

\subsubsection{Unique determination of the class of a divisor
by a representing $G$-form}
\label{subsubsection:UniqueDetermination}
 Let $G$ and $n$ be as above.
 Suppose that divisors $D$ and $D'$ of degree $\frac{1}{2}\deg G$
are represented by
$k$-equivalent $G$-forms $X$ and $X'$, respectively. We
claim that
$D$ and
$D'$ belong to the same divisor class. To prove the claim we may assume
without loss of generality that $X=X'$. Write
$X=uv=u'v'$
where $u$ and $u'$ are column vectors, $v$ and $v'$ are
row vectors, and the entries of
$u$ (resp.~$v$, $u'$, $v'$) form a $k$-basis for $L(D)$
(resp.~$L(G-D)$, $L(D')$, $L(G-D')$).
There
exists a unique nonzero meromorphic function $f$ on $\curve$ such that
$u'=fu$ and $v'=f^{-1}v$, hence 
$$L(D)=f^{-1}\cdot L(D')=L(D'+(f))\subset
L(\min(D,D'+(f))),$$
and hence $D=D'+(f)$ by
Lemma~\ref{Lemma:FundamentalUniquenessPrinciple}. The
claim is proved.

\subsubsection{Construction of a divisor represented by a given $G$-form}
\label{subsubsection:MeetsEvery}
 Let $G$ and $n$ be as above.
Let $X$ be a $G$-form. Choose
a column
$u$ and a row $v$ of $X$ each with $k$-linearly independent
entries.
Let $f$ be the entry common to
$u$ and
$v$; then $f$ does not vanish identically. 
Consider now the effective divisors
$$D=G+\min_j (v_j),\;\;\;E=G+\min_i (u_i),\;\;\;F=G+(f).$$
We claim that $X$ represents the divisor $D$.
In any case, since $X$ is of rank $\leq 1$ we have $X=uv/f$ and hence
$$D+E\geq F,\;\;\;\deg D+\deg E\geq \deg F=\deg G.$$
Since the entries of $v$ are $k$-linearly independent
and belong to $L(G-D)$, we have $\ell(G-D)\geq n$ 
and hence
$\frac{1}{2}\deg G\geq\deg D$ by
Lemma~\ref{Lemma:FundamentalUniquenessPrinciple} . Similarly we have
$\frac{1}{2}\deg G\geq\deg E$. It follows that
$$\deg D=\deg E=\frac{1}{2}\deg G=\frac{1}{2}\deg F,\;\;\;D+E=F.$$
In turn it follows that the entries of $v$ form a $k$-basis
for
$L(G-D)$ and that the entries of $u/f$ form a $k$-basis for
$L(G-E+(f))=L(D)$. Therefore $X=uv/f$ does indeed
represent the divisor
$D$. The claim is proved.

\begin{Proposition}\label{Proposition:Dictionary}
 Let $G$ be a divisor of even degree such that
$\frac{1}{2}\deg G\geq 2g$. There exists a unique
bijective correspondence 
$$\left\{\begin{array}{l}
\mbox{\textup{$k$-equivalence }}\\
\mbox{\textup{classes of $G$-forms}}
\end{array}\right\}\leftrightarrow\left\{\begin{array}{l}
\mbox{\textup{divisor classes}}\\
\mbox{\textup{of degree $\frac{1}{2}\deg G$}}
\end{array}\right\}
$$
with respect to which, for any  $G$-form
$X$ and any divisor
$D$ of degree $\frac{1}{2}\deg G$, the $k$-equivalence class of
$X$ corresponds to the divisor class of $D$
  if and only
if $X$ represents $D$.
\end{Proposition}
\proof All the hard work is done; we just have to
sum up. For each divisor
$D$ of degree
$\frac{1}{2}\deg G$ arbitrarily fix a column
(resp.~row) vector $u_D$ (resp.~$v_D$) with entries
forming a $k$-basis of
$L(D)$ (resp.~$L(G-D)$). The product $u_Dv_D$ is a
$G$-form the 
$k$-equivalence class of which depends only on $D$,
not on the choice of
$u_D$ and $v_D$ (\S\ref{subsubsection:DivMatrixRep}). The
map
$D\mapsto u_Dv_D$ sends divisor classes into \linebreak
$k$-equivalence classes
(\S\ref{subsubsection:DivMatrixRep}). The map $D\mapsto
u_Dv_D$ sends distinct divisor classes into distinct
$k$-equivalence classes
(\S\ref{subsubsection:UniqueDetermination}). The image
of the map \linebreak
$D\mapsto u_Dv_D$ meets every $k$-equivalence class of
$G$-forms (\S\ref{subsubsection:MeetsEvery}).
\qed

\subsubsection{Remark} 
The dictionary provided by
Proposition~\ref{Proposition:Dictionary} is essentially just a chapter
of  algebro-geometrical folklore, cf.\
\cite[Prop.~1.1]{EisenbudKohStillman}. 
The Eisenbud-Koh-Stillman paper greatly inspired us. Our
presentation of the dictionary is a much simplified  version of
the presentation in
\cite{Anderson}.

\subsection{Matrix representation of divisor class
addition and subtraction} 
\label{subsection:AdditionAndSubtraction}

\subsubsection{A theorem of Mumford}
By \cite[Thm.~6, p. 52]{Mumford} we have
\begin{equation}\label{equation:FirstMumfordTheorem}
\left.\begin{array}{rcl}
\deg D&\geq& 2g+1\\
\deg E&\geq&
2g
\end{array}\right\}\Rightarrow L(D+E)= \mbox{$k$-span
of $L(D)\cdot L(E)$}
\end{equation}
for all divisors $D$ and $E$.

\subsubsection{Kronecker products}
Given a
$p$ by
$q$ matrix
$A$ and an
$r$ by
$s$ matrix $B$ both with entries in some ring $R$, the
{\em Kronecker product}
$A\circ B$ is defined to be the $pr$ by $qs$ matrix with
entries in $R$ admitting a
decomposition into $r$ by $s$ blocks of the form
$$A\circ B=\left[\begin{array}{ccc}
&\vdots&\\
\dots&A_{ij}B&\dots\\
&\vdots
\end{array}\right].$$
The
Kronecker product of matrices is compatible with ordinary
matrix multiplication in the sense that
$$(A\circ B)(X\circ Y)=(AX)\circ (BY)$$
whenever
$AX$ and $BY$ are defined.

\begin{Proposition}\label{Proposition:DivisorAddition}
Let divisors $G$, $G'$, $D$ and $D'$ be given subject to the
following conditions:
$$
\deg
G=2\cdot\deg D,\;\;\;\deg G'=2\cdot\deg D',$$
$$\min\left(\frac{1}{2}\deg G,\frac{1}{2}\deg G'\right)\geq
2g,\;\;\;
\max\left(\frac{1}{2}\deg G,\frac{1}{2}\deg G'\right)\geq
2g+1.$$
Put
$$n=\frac{1}{2}\deg G-g+1,\;\;\;n'=\frac{1}{2}\deg G'-g+1,$$
and
$$n''=\frac{1}{2}\left(\deg G+\deg G'\right)-g+1=n+n'+g-1.$$
Fix a $G$-form $X$ representing $D$ and a $G'$-form $X'$
representing $D'$. Let $P$ and
$Q$ be any
$nn'$ by $nn'$
permutation matrices and consider the block decomposition
$$P\left(X\circ X'\right)Q=
\left[\begin{array}{cc}
a&b\\
c&d\end{array}\right]
$$
where the block $d$ is $n''$ by $n''$ and the other blocks
are of the appropriate sizes. (i) For some
$P$ and
$Q$ the corresponding block
$d$ is $k$-general. 
(ii) For any $P$ and $Q$ such that
$d$ is $k$-general, $d$ is a 
$(G+G')$-form representing
$D+D'$.  
\end{Proposition}
\proof Write
$X=uv$ and $X'=u'v'$
where $u$ (resp.~$v$, $u'$, $v'$) is a column (resp.~row,
column, row) vector with entries forming a $k$-basis of
$L(D)$ (resp.~$L(G-D)$, $L(D')$, $L(G'-D')$).

(i) By Mumford's theorem
(\ref{equation:FirstMumfordTheorem}) the entries of
$u\circ u'$ span
$L\left(D+D'\right)$ over $k$ and hence
for some permutation matrix $P$  the last $n''$ entries of
 the vector $P\left(u\circ u'\right)$ form
a
$k$-basis of
$L\left(D+D'\right)$. Similarly, the entries
of $v\circ v'$ span
$L\left(G+G'-D-D'\right)$ over
$k$ and for some permutation matrix $Q$ the last $n''$
entries of
$\left(v\circ v'\right)Q$ form a
$k$-basis of 
$L\left(G+G'-D-D'\right)$. 
With $P$
and $Q$ thus chosen the block $d$ is a
$(G+G')$-form  representing 
$D+D'$ and {\em a fortiori} $k$-general.

(ii) Now suppose we are given
$P$ and
$Q$ such that $d$ is $k$-general. Then
the last $n''$ entries of $P\left(u\circ u'\right)$ are
forced to be
$k$-linearly independent and hence to form a $k$-basis of
$L\left(D+D'\right)$. Similarly the last $n''$
entries $\left(v\circ v'\right)Q$ are
forced to form a
$k$-basis of
$L\left(G+G'-D-D'\right)$.
Then the block $d$ is indeed a
$\left(G+G'\right)$-form representing
$D+D'$. \qed

\begin{Lemma}\label{Lemma:Duality}
Let $E$ be a nonzero effective divisor.
Let $\RR_E$ be the ring consisting of the meromorphic functions on $C$
regular in a neighborhood of the support of $E$ and let $\II_E\subset
\RR_E$ be the ideal consisting of functions vanishing to order at least
$E$. Then there exists a $k$-linear functional
$$\sigma:\RR_E\rightarrow k$$ factoring through the quotient
$\RR_E/\II_E$ such that the induced $k$-bilinear
 map
\begin{equation}\label{equation:PerfectPairing}
\left((a\bmod{\II_E},b\bmod{\II_E})\mapsto
\sigma(ab)\right):\RR_E/\II_E\times \RR_E/\II_E\rightarrow
k
\end{equation} 
is a perfect pairing of $(\deg E)$-dimensional vector spaces
over
$k$.
\end{Lemma}
\proof Choose
any meromorphic differential
$\omega$ on
$C$ such that 
$$\ord_x \omega+\ord_x E=0$$
for all points $x\in \supp E$, where $\ord_x$ abbreviates  ``order of
vanishing at $x$''. Then, so we claim, the
$k$-linear functional
$$\left(a\mapsto \sum_{x\in \supp E}\Res_x(a\omega)
\right):\RR_E\rightarrow k$$
has all the desired properties. The proof of the claim is an exercise in
residue calculus we can safely omit.
\qed

\begin{Lemma}
\label{Lemma:CompressionFunctional}
Let $G$  and $E$ be divisors such that
$$\deg G\equiv 0\bmod{2},\;\;\;E>0,\;\;\;\frac{1}{2}\deg G-\deg
E>2g-2.$$ There exists a
$k$-linear functional 
$$\rho:L(G)\rightarrow k$$
factoring through the quotient $\frac{L(G)}{L(G-E)}$
such that for all divisors
$D$ of degree
$\frac{1}{2}\deg G$ 
the induced $k$-bilinear map
$$
\textstyle\left((a+L(D-E),b+L(G-D-E))\mapsto
\rho(ab)\right):
\frac{L(D)}{L(D-E)}\times
\frac{L(G-D)}{L(G-D-E)}
\rightarrow k
$$
is  a perfect pairing of $(\deg E)$-dimensional vector
spaces over $k$. 
\end{Lemma}
\proof 
For each divisor $D$ choose a meromorphic function
$f_D$ on $C$ such that 
$$\ord_x f=\ord_x D$$
for all points $x\in \supp E$;
then we have
$$L(D)\subseteq f^{-1}_D \RR_E,\;\;\;L(D)\cap
f^{-1}_D\II_E=L(D-E),$$ 
where $\RR_E$ and $\II_E$ are as defined in Lemma~\ref{Lemma:Duality}.
Now choose a
$k$-linear functional
$\sigma:\RR_E\rightarrow k$ such that the
pairing (\ref{equation:PerfectPairing}) is perfect and put
$$\rho=(x\mapsto \sigma(f_Gx)):L(G)\rightarrow k,$$
thereby defining a $k$-linear functional factoring
through the quotient $\frac{L(G)}{L(G-E)}$.  
Suppose now that $\deg D=\frac{1}{2}\deg G$
and consider the commutative diagram
$$\begin{array}{ccccccr}
\frac{L(D)}{L(D-E)}&\times&\frac{L(G-D)}{L(G-D-E)}&
\stackrel{\times}\rightarrow&\frac{L(G)}{L(G-E)}&
\stackrel{\rho}
\rightarrow&
k\\
\downarrow&&\downarrow &&\downarrow&&
\downarrow
\\
\RR_E/\II_E&\times &\RR_E/\II_E
&\stackrel{\times}\rightarrow
&\RR_E/\II_E&\stackrel{\sigma}\rightarrow& k
\end{array}
$$
where the vertical arrows are induced by multiplication by
$$f_D,\;\;f_G f_{D}^{-1},\;\;f_G,\;\;1,$$
respectively. By construction all the vertical arrows are
injective and of course the last is  bijective. Further,
the source of each vertical arrow other than the
last is of dimension over $k$ equal to $\deg E$ by Riemann-Roch.
Therefore
all the vertical arrows are bijective and hence $\rho$ has the
desired nondegeneracy property.
\qed

\subsubsection{Compression functionals}
In the situation of Lemma~\ref{Lemma:CompressionFunctional}
we call 
$\rho:L(G)\rightarrow k$ an {\em
$E$-compression functional}. 
To make the calculations below run smoothly it is convenient to introduce
in this context the following notation. Given any matrix $Y$ with entries
in
$L(G)$, let $\rho Y$ be the result of applying $\rho$ entrywise to $Y$,
i.~e.,
the matrix with entries in
$k$ defined by the rule $(\rho Y)_{ij}=\rho Y_{ij}$.

\begin{Proposition}\label{Proposition:DivisorSubtraction}
Let $G$ and $E$ be divisors such that
$$\deg G\equiv 0\bmod{2},\;\;\;E>0,\;\;\;
\frac{1}{2}\deg G-\deg E\geq 2g$$
and put 
$$n=\frac{1}{2}\deg G-\deg
E-g+1,\;\;\;
n'=\frac{1}{2}\deg G-g+1=n+\deg E.$$
 Let 
$\rho:L(G)\rightarrow k$ be an
$E$-compression functional. Let
$D$ be a divisor of degree $\frac{1}{2}\deg
G$ and let $X$ be a $G$-form representing
$D$.  Let
$P$ and
$Q$ be any
$n'$ by
$n'$ permutation matrices and consider the block
decomposition
$$PXQ=\left[\begin{array}{cc}a&b\\
c&d\end{array}\right]$$
where the block $a$
is $\deg E$ by $\deg E$, the block $d$ is $n$ by
$n$ and the other blocks are of the appropriate
sizes.  
(i) For some
$P$ and
$Q$ we have $\det \rho a\neq 0$. (ii) For any 
$P$ and
$Q$ such that $\det \rho a\neq 0$
the matrix $z$ defined by
the rule
$$\left[\begin{array}{cc} w &x\\
y&z\end{array}\right]=\left[\begin{array}{cc}
1&0\\ -(\rho c)(\rho
a)^{-1}&1\end{array}\right]\left[\begin{array}{cc} a&b\\
c&d\end{array}\right]
\left[\begin{array}{cc}
1 &-(\rho a)^{-1}(\rho b)\\
0&1\end{array}\right]$$
is
a
$(G-2E)$-form representing
$D-E$. 
\end{Proposition}
\proof 
Write $X=uv$
where $u$ (resp.~$v$) is a column (resp.~row) vector with
entries forming a
$k$-basis of
$L(D)$ (resp.~
$L(G-D)$).

(i)   For suitably chosen permutation matrices $P$
and $Q$ the first
$\deg E$ entries of
the column vector $Pu$ (resp.~row vector $vQ$) project
to a
$k$-basis of the quotient
$\frac{L(D)}{L(D-E)}$ (resp.~$\frac{L(G-D)}{L(G-D-E)}$).
For such $P$ and $Q$ we have $\det
\rho a\neq 0$ by
definition of
$E$-compression
functional.

(ii) After replacing $X$ by $PXQ$ we may assume that
$P=Q=1$. After replacing $X$ by the $k$-equivalent
matrix
$\left[\begin{array}{cc} w&x\\
y&z\end{array}\right]$ we may assume that
$$\det \rho a\neq 0,\;\;\;\rho b=0,\;\;\;\rho c=0,$$
in which case our task is simply to show that the block $d$ is a
$(G-2E)$-form. By definition of $E$-compression functional
the first
$\deg E$ entries of
$u$ (resp.~$v$) project to a $k$-basis of the quotient
$\frac{L(D)}{L(D-E)}$ (resp.~$\frac{L(G-D)}{L(G-D-E)}$).
Also by definition of  $E$-compression functional the last
$n$ entries of $u$ must belong to $L(D-E)$ and since
$k$-linearly independent must form a
$k$-basis of
$L(D-E)$.
Similarly the last $n$ entries of $v$ must form a
$k$-basis of $L(G-D-E)$. Therefore $d$ is indeed a
$(G-2E)$-form representing 
$D-E$.
\qed

\subsection{Completion of the construction}

\subsubsection{Candidate for the Jacobian}
\label{subsubsection:JacobianCandidate}
 Fix an effective divisor $E$ of degree $\geq 2g+1$ and put
$$S=C\setminus \supp E,\;\;\;A=H^0(S,\OO_C),\;\;\;
n=\ell(E)=\deg E-g+1,\;\;\;
L=L(2E).
$$
The  projective algebraic
variety $J$ of
$k$-proportionality classes of Jacobi matrices of type
$\left(k,n,A,L\right)$ is our candidate for the Jacobian of
$C$.

\subsubsection{Candidate for the Abel map} 
\label{subsubsection:AbelCandidate}
For each divisor $D$ of degree zero
arbitrarily fix a
$2E$-form
$X_D$ representing
$D+E$. Now a $2E$-form is the same thing as a Segre matrix
of type $(k,n,A,L)$. By
Proposition~\ref{Proposition:Dictionary} it follows that
the map $D\mapsto X_D$ puts the divisor classes of degree zero
in bijective correspondence with the
$k$-equivalence classes of Segre matrices of type
$(k,n,A,L)$. For
each divisor
$D$ of degree zero let
$Z_D$ be the image of $X_D$ under the abstract Abel map.
 By Theorem~\ref{Theorem:ProjectiveParameterization} it
follows that the map
$D\mapsto Z_D$ puts the classes of divisors of degree zero
into bijective correspondence with the points of
$J$. The bijective map from classes of divisors of degree
zero to 
$J$ induced by the map
$D\mapsto Z_D$ is our candidate for the Abel map. 

\begin{Lemma}\label{Lemma:InversionAlgebraicity}
For all
divisors $D$ of degree zero,
$X_{-D}$ is $k$-equivalent to $X_D^T$,
and (hence) $Z_{-D}$ is $k$-proportional to $Z_D^T$.
\end{Lemma}
\proof This boils down to the
transpose symmetry (\ref{equation:AbeliantTransposeSymmetry})
of the abeliant.
\qed

\begin{Lemma}\label{Lemma:RoughGroupLaw}
Fix an
$E$-compression functional
$\rho:L(4E)\rightarrow k$. 
 Fix Segre matrices $X$ and $X'$ of type $(k,n,A,L)$.
Fix a divisor $D$ (resp.~$D'$) such that $X$ (resp.~$X'$)
is $k$-equivalent
to $X_D$ (resp.~$X_{D'}$). Let
$P$ and
$Q$ be any
$n^2$ by
$n^2$ permutation matrices
 and consider the block
decomposition
$$P\left(X\circ X'\right)Q
=\left[\begin{array}{ccc}
\bullet &\bullet&\bullet\\
\bullet&a&b\\
\bullet&c&d
\end{array}\right]
$$
where 
the block $a$ is $\deg E$ by $\deg E$, the block $d$ is
$n$ by
$n$, the other blocks are of the appropriate
sizes, and the bullets hold places for blocks the contents
of which do not concern us.   Further, consider the
block-decomposed matrix
$$\left[\begin{array}{cc} w&x\\
y&z\end{array}\right]= \det \rho a\cdot
 \left[\begin{array}{cc} \det \rho a&0\\
-(\rho c)(\rho
a)^\star&\det
\rho a\end{array}\right]\left[\begin{array}{cc}
 a&b\\ c&d\end{array}\right]
\left[\begin{array}{cc}
\det \rho a &-(\rho a)^\star(\rho b)\\
0&\det \rho a\end{array}\right].$$
(i) For some
$P$ and
$Q$ the corresponding block $z$ is $k$-general. (ii)
For any $P$ and $Q$ such that the corresponding block
 $z$ is $k$-general, 
$z$ is a Segre matrix of type $(k,n,A,L)$ and moreover $z$ is
$k$-equivalent to $X_{D+D'}$. 
\end{Lemma}
\proof (i) By
Propositions~\ref{Proposition:DivisorAddition} and
\ref{Proposition:DivisorSubtraction} there exist
$P$ and
$Q$ such that the following hold:\\
\begin{itemize}
\item 
$\left[\begin{array}{cc} a&b\\
c&d\end{array}\right]$ is a $4E$-form representing
$D+D'+2E$.\\
\item $\det \rho a\neq 0$.\\
\item $z$ is a
$2E$-form representing $D+D'+E$.\\
\end{itemize}
{\em A fortiori} $z$ is $k$-general. \\

(ii) By hypothesis we have
$$\det\rho w\neq 0,\;\;\;\rho
x=0,\;\;\;\rho y=0,$$
and there exists a factorization
$$\left[\begin{array}{cc} w&x\\
y&z\end{array}\right]=\left[\begin{array}{c}
p\\q\end{array}\right]\left[\begin{array}{cc}
r&s\end{array}\right]
$$
where the entries of the column vector (resp.~row vector)
on the right belong to
$L(D+D'+2E)$ 
(resp.~$L(-D-D'+2E)$), the blocks $p$ and $r$ are
vectors of length
$\deg E$, and the blocks $q$ and $s$ are vectors of length
$n$. 
By definition of
$E$-compression
functional (Lemma~\ref{Lemma:CompressionFunctional}) it
follows that the entries of $p$ (resp.~$r$) project to a
$k$-basis of the quotient
$\frac{L(D+D'+2E)}{L(D+D'+E)}$ 
(resp.~
$\frac{L(-D-D'+2E)}{L(-D-D'+E)}$).
Also by definition of
$2E$-compression
functional 
it follows that the entries of $q$ (resp.~$s$)
belong to
$L(D+D'+E)$ (resp.~$L(-D-D'+E)$). Finally,
since
$z=qs$ is $k$-general, the entries of $q$ (resp.~$s$)
must be
$k$-linearly independent, and hence the entries of $q$
(resp.~$s$) must form a $k$-basis of $L(D+D'+E)$
(resp.~$L(-D-D'+E)$). Therefore the block
$z$ is indeed a $2E$-form representing 
$D+D'+E$ and hence $k$-equivalent to $X_{D+D'}$.
\qed

\begin{Lemma}\label{Lemma:AdditionAlgebraicity}
Fix an
$E$-compression functional
$\rho:L(4E)\rightarrow k$. Fix Jacobi matrices
$Z$ and $Z'$ of type $(k,n,A,L)$ with discriminants
$\Delta$ and $\Delta'$, respectively. Fix a divisor $D$
(resp.~$D'$) of degree zero such that $Z$ (resp.~$Z'$) is
$k$-proportional to $Z_D$ (resp.~$Z_{D'}$).  For any
$\sbold,\sbold'\in S^{\{1,\dots,n+1\}}$ and any
$n^2$ by
$n^2$ permutation matrices $P$ and $Q$ consider the block
decomposition
$$P\left(\left(Z\Vert_\sbold\right)\circ
\left(Z'\Vert_{\sbold'}\right)\right)Q
=\left[\begin{array}{ccc}
\bullet&\bullet&\bullet\\
\bullet&a&b\\
\bullet&c&d
\end{array}\right],
$$
where 
the block $a$ is $\deg E$ by $\deg E$, the block $d$ is
$n$ by
$n$, the other blocks are of the appropriate
sizes, and the bullets hold places for blocks the contents
of which do not concern us.  Consider the block
decomposed matrix
$$\left[\begin{array}{cc} w&x\\
y&z\end{array}\right]=\Delta\Vert_\sbold\cdot
\Delta'\Vert_{\sbold'}\cdot \det \rho
a\cdot \left[\begin{array}{cc} \det \rho a&0\\
-(\rho c)(\rho
a)^\star&\det
\rho a\end{array}\right]\left[\begin{array}{cc}
 a&b\\ c&d\end{array}\right]
\left[\begin{array}{cc}
\det \rho a &-(\rho a)^\star(\rho b)\\
0&\det \rho a\end{array}\right]$$
and finally put
$$Z''
=\abel_{\ell=0}^{n+1}z^{(\ell)}.$$
(i) There exist
$\sbold$,
$\sbold'$,
$P$ and
$Q$ such that
$Z''$ does not vanish identically. (ii) For any $\sbold$,
$\sbold'$,
$P$ and
$Q$ such that $Z''$ does not vanish identically,
 $Z''$ is a Jacobi matrix of type
$(k,n,A,L)$ and moreover $Z''$ is $k$-proportional to
$Z_{D+D'}$.
\end{Lemma}
\proof
(i) By Proposition~\ref{Proposition:Clincher} there
exists $\sbold\in S^{\{1,\dots,n+1\}}$
(resp.~$\sbold'\in S^{\{1,\dots,n+1\}}$) such that
$\Delta\Vert_\sbold$
(resp.~$\Delta'\Vert_{\sbold'}$) is a nonzero scalar
and hence the corresponding partially specialized matrix
$Z\Vert_\sbold$ (resp.~$Z'\Vert_{\sbold'}$)
is $k$-equivalent to $X_D$ (resp.~$X_{D'}$).
By Lemma~\ref{Lemma:RoughGroupLaw} there exist $P$ and $Q$
such that the block
$z$ is $k$-general. Finally, $Z''$  does not vanish
by Proposition~\ref{Proposition:GeneralTypeTilde}.

(ii) By 
Proposition~\ref{Proposition:Clincher} and hypothesis
the partial specialization
$\Delta\Vert_\sbold$ \linebreak 
(resp.~$\Delta'\Vert_{\sbold'}$) is a nonzero scalar
and hence the
 corresponding partial specialization
$Z\Vert_\sbold$ (resp.~$Z'\Vert_{\sbold'}$) is
a Segre matrices of type $(k,n,A,L)$ that is
$k$-equivalent to $X_D$ (resp.~$X_{D'}$).
By hypothesis $\det \rho a$ is a
nonzero scalar and hence $z$ is by 
Lemma~\ref{Lemma:RoughGroupLaw}  a Segre matrix of
type $(k,n,A,L)$ that is $k$-equivalent to
$X_{D+D'}$. Finally, since
$Z''$ is the image of $z$ under the abstract Abel map, 
$Z''$ is $k$-proportional to $Z_{D+D'}$.
\qed

\begin{Theorem}\label{Theorem:Main}
 There
is exactly one way to equip our candidate for the
Jacobian (\textup{\S\ref{subsubsection:JacobianCandidate}}) with
the structure of algebraic group so that our candidate
for the Abel map (\textup{\S\ref{subsubsection:AbelCandidate}})
becomes a group homomorphism. (Thus our candidates become {\em the}
Jacobian and {\em the} Abel map.)
\end{Theorem}
\proof Since our candidate for the
Abel map is bijective, the set underlying $J$ comes
canonically equipped with a group law. The only issue
remaining to be resolved is whether or not  that  group law
is algebraic, i.~e., expressible Zariski-locally by regular
functions. Well, by Lemma~\ref{Lemma:InversionAlgebraicity}
the inversion operation in
$J$ is algebraic, and
by Lemma~\ref{Lemma:AdditionAlgebraicity}
the addition operation in $J$ is algebraic.
We're done.
\qed

\subsection{Remark}\label{Remark:Last}
To give some indication of how the complexity
of our construction of
$J$ grows as a function of the genus of
$C$, we make the following observation.
Suppose that $C$ is a nonsingular plane algebraic curve of
degree
$d\geq 3$ and hence genus $\frac{(d-1)(d-2)}{2}>0$ with
defining equation
$F=F(x,y,z)\in k[x,y,z]$.
By the method of proof of
Theorem~\ref{Theorem:Main} the divisor classes of $C$ of degree zero
can be put in natural bijective correspondence with the
$k$-equivalence classes of Segre matrices of type
\begin{equation}\label{equation:PlaneType}
\left(k,\frac{d(d-1)}{2},k[x,y,z]/(F),\{\mbox{forms of
degree
$2d-4$}\}/(F)\right)
\end{equation}
and in turn the Jacobian of $C$ can be identified with the projective
variety of
$k$-proportionality classes of Jacobi matrices of type
(\ref{equation:PlaneType}). 

\section{Acknowledgements}
I thank Dinesh
Thakur for comments on preliminary  drafts
of this paper. I thank Joel Roberts for
discussions concerning ideals generated by two by two
minors.   I thank Hendryk Lenstra for
comments on a preliminary draft of this paper, and in particular for
pointing out Pila's paper to me. I thank Jeremy Teitelbaum for a
conversation helpful for devising the example of
\S\ref{Remark:Last}. I
thank the referee for constructive criticism.

\end{document}